\documentclass[11pt,letterpaper]{amsart}

\usepackage{tikz}
\usetikzlibrary{arrows.meta}
\usepackage{url}
\usepackage{graphicx}
\usepackage{amssymb}
\usepackage[T1]{fontenc}
\usepackage{mathscinet}
\usepackage{enumitem}

\allowdisplaybreaks%

{\theoremstyle{plain}%
 \newtheorem{theorem}{Theorem}
 
 \newtheorem{lemma}{Lemma}%
}
{\theoremstyle{remark}
\newtheorem{remark}{Remark}
}
{\theoremstyle{definition}

}

\title[On a problem due to Erd\H{o}s concerning $w_1(n)$]{On a problem due to Erd\H{o}s concerning $w_1(n)$}
\author{John M.\ Campbell}
\address{Department of Mathematics and Statistics, Dalhousie University,
Halifax, NS B3H 4R2, Canada}
\email{jh241966@dal.ca}
 


\begin{document}

\begin{abstract}
 Adopting notation from Erd\H{o}s [\emph{Congress.\ Numer.}, 1980], we define $w_1(n) = w(n)$ for positive integers $n$ so that $w(1) = 
 0$ and so that: If $n > 1$ and $n = \prod p^{\alpha_p}$, then $w(n) = \sum_{p^{\alpha_p} \parallel n} \frac{n}{p^{\alpha_p}}$. By then 
 defining $V(x)=\bigl|\{w(n):1\le n\le x\}\bigr|$ for $x \geq 1$, Erd\H{o}s considered whether or not $V(x) = o(x)$ and left this as an 
 open problem. This same problem was also given as an open problem in Guy's famous textbook on unsolved problems in number theory 
 and seems to have remained open. We solve this problem by proving that $V(x) \neq o(x)$, and, moreover, by proving the stronger result 
 such that $\liminf_{x \to \infty} \frac{V(x)}{x} > 0$. 
 \end{abstract}

\maketitle

\section{Introduction and background}
 Define the integer sequence $(w(n) : n \in \mathbb{N})$ so that 
\begin{equation}\label{definewn} 
 w(n) = \begin{cases} 
 \sum_{p^{\alpha_p}\parallel n} \frac{n}{p^{\alpha_p}} & \text{if $n > 1$ and $n = \prod p^{\alpha_p}$; and} \\
 0 & \text{if $n = 1$,} 
 \end{cases}
\end{equation}
 and define $V(x)=\bigl|\{w(n):1\le n\le x\}\bigr|$. For $n > 1$, the function in \eqref{definewn} is the function denoted by $w_{1}(n)$ by 
 Erd\H{o}s \cite{Erdos1980}, who also considered its iterates. Such iterates will not be considered in this paper, so we write $w(n) = 
 w_1(n)$. We also adopt the empty-sum 
 convention such that $w(1) = 0$. A problem introduced by Erd\H{o}s in 1980 \cite{Erdos1980} is given by whether or not $V(x) = o(x)$, and 
 Erd\H{o}s described this problem by stating it ``perhaps is not hopeless.'' This same problem is also 
 given as an unsolved problem in Guy's 
 famous text on unsolved problems in number theory \cite[p.\ 64]{Guy1994}, and it seems that this same problem has remained open, 
 in view of the extant references citing Erd\H{o}s's work \cite{Erdos1980}, and in view of the extant material in the entries in the On-line 
 Encyclopedia of Integer Sequences\footnote{See \url{https://oeis.org/}.} that agree (up to initial terms) with the integer sequence
 defined in \eqref{definewn}, namely, the OEIS entries {\tt A028236} and {\tt A066504}. We succeed, based on our extensive interactions 
 with GPT‑6 Astra Pro, in solving Erd\H{o}s's problem, by proving that $V(x) \neq o(x)$. Furthermore, we prove the stronger result such that 
 $\liminf_{x \to \infty} \frac{V(x)}{x} > 0$. 

\subsection{Brief outline}
Our proof uses integers of the form $n=etrp\leq x$, where $e\equiv2\pmod4$ and $t$, $r$, and $p$ are primes in suitably separated ranges. 
 We construct a set of such integers whose cardinality is bounded below by a fixed positive multiple of $x$, and we bound the number of 
 ordered pairs of its elements having the same value under $w$ by a fixed multiple of $x$. The desired lower bound for $V(x)$ then 
 follows from the Cauchy--Schwarz inequality. Writing $a=etr$, the identity $w(ap)=p\,w(a)+a$ reduces the equality $w(ap)=w(b 
 q)$ to a linear equation in the primes $p$ and $q$, to which we
apply a uniform sieve estimate.

 The two smaller prime factors have different roles in estimating the resulting weights. For fixed distinct $c$ and $d$ in this 
 construction, writing $a=cr$ and $b=ds$, every positive common divisor of $w(cr)$, $w(ds)$, and $cr-ds$ divides the nonzero integer 
 $c^2w(d)-d^2w(c)$. The chosen ranges ensure that this integer is smaller in absolute value than the primes $r$ and $s$, allowing an 
 average over $r$ and $s$ through the Brun--Titchmarsh inequality. Writing $c=et$ and $d=fu$, we then average over $t$: For fixed 
 $e$, $f$, and $u$, the same integer is the value of a quadratic polynomial in $t$ whose
content is $\gcd(e,f)^2$. A uniform estimate for polynomial values at primes, followed by a harmonic gcd estimate, controls 
 the remaining weight. Related polynomial averages control
 the prime-factor weights contributed by the sieve. The case $c=d$ is treated separately.

\section{Preliminaries}\label{secprelim}
 Let $D\subseteq[0,\infty)$ be unbounded, and let $F,G\colon D \to\mathbb R$ be functions that may depend on parameters $\rho_1, 
 \rho_2, \ldots, \rho_{n}$. Assume that: For each fixed admissible choice of these parameters, the relation $G(x) \geq 0$ holds for all 
 sufficiently large $x \in D$. We then write $F \ll_{\rho_1, \rho_2, \ldots, \rho_n} G$ as $x \to \infty$ through $D$ to indicate that: For each 
 fixed and admissible choice of $\rho_{1}$, $\rho_{2}$, $\ldots$, $\rho_{n}$, there are constants $C = C(\rho_1, \rho_2, \ldots, \rho_n) > 
 0$ and $x_0 = x_0(\rho_1, \rho_2, \ldots, \rho_n)\geq 0$ (possibly depending on the given parameters) such that $|F(x)| \leq C \, 
 G(x)$ holds for all $x \in D$ such that $x \geq x_0$. Observe that $C$ and $x_0$ are independent of $x$. Also observe that no 
 uniformity with respect to the given parameters is asserted. We may then define $\ll$ accordingly, via an empty set of parameters. 

 We require the use of a number of familiar and fundamental arithmetic functions. In this direction, we let $\omega(n)$ denote the
 number of distinct prime divisors of $n \geq 2$, with the convention such that $\omega(1) = 0$. We recall the 
 usual notation $\pi(x)$ for the prime-counting function giving the number of primes not exceeding $x$, 
 and we recall one of the most fundamental results in all of number theory, 
 namely, the \emph{prime number theorem} asserting that $\pi(x) \sim \frac{x}{\log x}$. 
 We also make use of the usual notation for the counting function for primes in arithmetic progressions, writing 
\begin{equation}\label{countarithmetic}
 \pi(x;d, a) = \big| \{ p : 
 p \ \text{is prime}, p \leq x, p \equiv a \bmod{d} \} \big| 
\end{equation}
 for $x \in \mathbb{R}_{\geq 0}$, $d \in \mathbb{Z}_{\geq 1}$, and $a \in \mathbb{Z}$. 
 Recalling that \emph{Euler's totient function} is defined so that 
 $
 \varphi(n) = \big| \{ a \in \mathbb{Z} : 1 \leq a \leq n, \gcd(a, n) = 1 \} \big|$, 
 a well known generalization of the prime number
 theorem is such that $\pi(x$; $d$, $a)$ $ \sim$ $ \frac{x}{\varphi(d) \log x}$ (i.e., for fixed $d$ and $a$ 
 satisfying $\gcd(a,d)=1$ and as $x \to \infty$), 
 but we require an explicit bound for $ \pi(x;d,a)$ given below. 

 A fundamental property concerning the arithmetic function defined in \eqref{countarithmetic} is the \emph{Brun--Titchmarsh inequality}. 
 A formulation of this result attributed to Montgomery and Vaughan \cite{MontgomeryVaughan1973} gives us that: If $d$ and $a$ are 
 positive integers that are relatively prime, then
\begin{equation}\label{BrunTitchmarsh}
 \pi(x;d,a) < \frac{2x}{\varphi(d) \log\left( \frac{x}{d} \right)} 
\end{equation}
 holds for all $x > d$, referring to Crandall and Pomerance's text for background \cite[p.\ 43]{CrandallPomerance2005}. We require the 
 relation in \eqref{BrunTitchmarsh} for our main construction. 

 We henceforth adopt the convention whereby variables denoted with $p$  in indexing sets for sums and products are prime,  and we 
 recall \emph{Euler's  product formula} for the totient function giving that 
\begin{equation}\label{Eulerprod}
 \varphi(n) = n \prod_{p \mid n} \left( 1 - \frac{1}{p} \right) 
\end{equation}
 for each positive integer $n$ (with the empty product being interpreted as being equal to $1$). 
 We also recall \emph{B\'{e}zout's identity} giving us that: For integers $A$ and $B$ that are not both equal to $0$, there exist integers 
 $u$ and $v$ such that $ A u + B v = \gcd(A, B)$. 

 We require some fundamental inequalities often applied in analysis and number theory, including 
 a formulation of \emph{H\"{o}lder's inequality} giving us that:
 For real numbers $p, q > 1$ such that 
 $\frac{1}{p} + \frac{1}{q} = 1$, 
 and for a positive integer $N$, 
 and for complex numbers $a_{1}$, $a_{2}$, $\ldots$, $a_{N}$
 and $b_{1}$, $b_{2}$, $\ldots$, $b_{N}$, 
 the inequality 
\begin{equation}\label{generalHolder} 
 \sum_{i=1}^{N} |a_{i} b_{i}| \leq 
 \left( \sum_{i=1}^{N} |a_{i}|^{p} \right)^{\frac{1}{p}} \left( \sum_{i=1}^{N} |b_{i}|^{q} \right)^{\frac{1}{q}} 
\end{equation}
 holds. 

 We also make use of some familiar terminology arising in the number-theoretic study of polynomials. For a nonzero polynomial $P$ with 
 integer coefficients, the \emph{content} of $P$, denoted by $\operatorname{cont}(P)$, 
 refers to the positive greatest common divisor of its coefficients. 
 A \emph{primitive polynomial} is a nonzero polynomial with integer coefficients whose coefficients
 have greatest common divisor $1$. 

\section{Main construction}
 We begin by briefly establishing some basic and required properties concerning the integer sequence in \eqref{definewn}. 

\begin{lemma}\label{gcdmwm}
 For each positive integer $n$, the equality $ \gcd(n, w(n)) = 1 $ holds. 
 \end{lemma}

\begin{proof} 
 This is immediate for $n = 1$. If $p^{\alpha} \parallel n > 1$, then all terms in the expansion of $w(n)$ displayed in \eqref{definewn}
 apart from $\frac{n}{p^{\alpha}}$ are divisible by $p$, whereas the exceptional term is not. 
 Thus no prime divisor of $n$ divides $w(n)$, proving the claim.
\end{proof}

\begin{lemma}\label{lemmawap}
 Suppose that $a$ is a positive integer and $p$ is a prime such that $p \nmid a$. Then the relation $ w(ap) = p \, w(a) + a $ holds. 
\end{lemma}

\begin{proof}
 This is a special case of the relation such that: If $a$ and $b$ are positive and relatively prime integers, then $ w(a b) = a \, w(b) + 
 b \, w(a)$, and, in turn, this follows in a direct way from the definition in \eqref{definewn}. 
\end{proof} 

\begin{lemma}\label{lemboundw}
 The upper bound such that $ w(n) \leq n \frac{\log n}{\log 2} $
 holds for each positive integer $n$. Moreover, if $2 \parallel n$, 
 then $\frac{n}{2} \leq w(n)$. 
\end{lemma}

\begin{proof}
 For $n=1$, the upper bound holds because $w(1)=0$ and $\log 1=0$. Suppose that $n>1$. The sum defining $w(n)$
has $\omega(n)$ terms, each at most $n$, so that 
$w(n)\leq n\,\omega(n)$.

 The product of the distinct prime divisors of $n$ divides $n$, and each of these primes is at least $2$. Consequently, 
 the relation $n\geq 2^{\omega(n)}$ holds. Taking logarithms gives $\omega(n)\leq \frac{\log n}{\log 2}$, and 
 therefore $ w(n)\leq n\,\omega(n)\leq\frac{n 
 \log n}{\log 2}$. 

 If $2\parallel n$, then the sum defining $w(n)$ contains
the term $\frac{n}{2}$, and every remaining term is nonnegative, 
 and hence the desired relation whereby $w(n) \geq \frac{n}{2}$.
\end{proof}

\subsection{The $\mathcal{B}$-function}
 Let $N$ denote an integer exceeding $1$. Let $\beta > 0$ and $K > 0$ be real. A key to our construction is given by the function 
 $ \mathcal{B}_{\beta,K}(N)$ defined so that 
\begin{equation}\label{definecalB}
 \mathcal{B}_{\beta,K}(N)=\prod_{p \mid N}(1 + K p^{-\beta}) 
\end{equation}
 and so that $ \mathcal{B}_{\beta,K}(1)=1$. 

\begin{lemma}\label{lemcalBineq}
 The relation $ \mathcal{B}_{\beta, K}(MN) 
 \leq \mathcal{B}_{\beta, K}(M) \mathcal{B}_{\beta, K}(N) $ holds for 
 integers $M, N \geq 1$. 
\end{lemma}

\begin{proof}
 It is immediate from \eqref{definecalB} that 
\begin{equation}\label{calBgcd}
 \mathcal{B}_{\beta, K}(M) \mathcal{B}_{\beta, K}(N)
 = \mathcal{B}_{\beta, K}(M N) \mathcal{B}_{\beta, K}(\gcd(M, N)) 
\end{equation}
 for integers $M, N \geq 1$, and from \eqref{calBgcd}, we obtain the desired relation. 
\end{proof}

\begin{lemma}\label{lemreciprocal}
 Let $\alpha$, $\delta$, and $\gamma$ denote fixed real numbers such that 
\begin{equation}\label{comparegamma}
 0 < \gamma < \alpha < \delta. 
\end{equation} 
 Then there exists a positive constant $C_{\alpha, \delta, \gamma}$
 depending only on $\alpha$, $\delta$, and $\gamma$
 such that: For each real number
\begin{equation}\label{Ygeq2}
 Y \geq 2, 
\end{equation}
 for each positive integer $h$ such that 
\begin{equation}\label{boundh}
 1 \leq h \leq Y^{\gamma} 
\end{equation} 
 and for each positive integer $b$ that is relatively prime with $h$, we have that 
\begin{equation}\label{explicitreciprocal}
 \sum_{\substack{ Y^{\alpha} < p \leq Y^{\delta} \\ p \equiv b \bmod {h} }} 
 \frac{1}{p} \leq \frac{C_{\alpha, \delta, \gamma}}{\varphi(h)}, 
\end{equation}
 so that 
 $$ \sum_{\substack{ Y^{\alpha} < p \leq Y^{\delta} \\ p \equiv b \bmod {h} }} 
 \frac{1}{p} \ll_{\alpha, \delta, \gamma} \frac{1}{\varphi(h)}, $$
 uniformly over all $Y$, $h$, and $b$ satisfying the given conditions. Moreover, 
 one may take 
 $$ C_{\alpha, \delta, \gamma} 
 = \frac{2}{(\delta - \gamma) \log 2} + 2 \log\left( \frac{\delta - \gamma}{\alpha - \gamma} \right). $$
\end{lemma}

\begin{proof}
 Since $Y > 1$ and $\gamma < \alpha$, we find that $h \leq Y^{\gamma} < Y^{\alpha}$. So, a given element $t$ in $[Y^{\alpha}, 
 Y^{\delta}]$ is such that $t > h$. 
 From the above formulation of the Brun--Titchmarsh inequality, 
 we find that 
\begin{equation}\label{BTconsequence}
 \pi(t;h,b) \leq \frac{2t}{ \varphi(h) \log\left( \frac{t}{h} \right) } 
\end{equation}
 for all $t \in [Y^{\alpha}, Y^{\delta}]$. 
 For each prime $p$ occurring in the sum shown in 
 \eqref{explicitreciprocal}, we find that 
\begin{equation}\label{sumbothsides} 
 \frac{1}{p} = \frac{1}{Y^{\delta}} + \int_{p}^{Y^{\delta}} \frac{dt}{t^2}. 
\end{equation}
 Summing both sides of \eqref{sumbothsides} and reversing the order of 
 (finite) summation and integration, we find that 
 $$ \sum_{\substack{ Y^{\alpha} < p \leq Y^{\delta} \\ p \equiv b \bmod {h} }} 
 \frac{1}{p} = 
 \frac{ \pi(Y^{\delta};h,b) - 
 \pi(Y^{\alpha};h,b) }{Y^{\delta}} + \int_{Y^{\alpha}}^{Y^{\delta}} \frac{ 
 \pi(t;h,b) - \pi(Y^{\alpha};h,b) }{t^2} \, dt, $$ and this, in turn, produces the equality 
\begin{equation}\label{simplifyint} 
 \sum_{\substack{ Y^{\alpha} < p \leq Y^{\delta} \\ p \equiv b \bmod {h} }} 
 \frac{1}{p} = 
 \frac{ \pi(Y^{\delta};h,b)
 }{Y^{\delta}} - \frac{ 
 \pi(Y^{\alpha};h,b) 
 }{Y^{\alpha}} + \int_{Y^{\alpha}}^{Y^{\delta}} \frac{ \pi(t;h,b) }{t^2} \, dt. 
\end{equation}
 From \eqref{simplifyint} and from the consequence in \eqref{BTconsequence} of the Brun--Titchmarsh inequality, we find that 
\begin{align*}
 \sum_{\substack{ Y^{\alpha} < p \leq Y^{\delta} \\ p \equiv b \bmod {h} }} 
 \frac{1}{p} & 
 \leq \frac{\pi(Y^{\delta}; h,b)}{Y^{\delta}} + \int_{Y^{\alpha}}^{Y^{\delta}} \frac{\pi(t;h,b)}{t^2} \, dt \\ 
 & \leq \frac{2}{ \varphi(h) \log\left( \frac{Y^{\delta}}{h} \right)} 
 + \frac{2}{\varphi(h)} \int_{Y^{\alpha}}^{Y^{\delta}} \frac{dt}{t \log\left( \frac{t}{h} \right)} \\ 
 & = \frac{2}{\varphi(h)} 
 \left( \frac{1}{\log\left( \frac{Y^{\delta}}{h} \right)} + \log\left( \frac{\log\left( \frac{Y^{\delta}}{h} \right)}{\log\left( 
 \frac{Y^{\alpha}}{h} \right)} \right) \right), 
\end{align*}
 noting that $\log\left( \frac{Y^{\alpha}}{h} \right) > 0$
 follows from \eqref{comparegamma} and 
\begin{equation}\label{loghleq} 
 \log h \leq \gamma \log Y, 
\end{equation}
 with \eqref{loghleq} being implied by \eqref{boundh}. Also, we find that 
\begin{align}
\begin{split}
 \frac{ \log\left( \frac{Y^{\delta}}{h} \right) }{ 
 \log\left( \frac{Y^{\alpha}}{h} \right) } 
 & = \frac{ \delta \log Y - \log h }{ \alpha \log Y - \log h } \\ 
 & = 1 + \frac{ (\delta - \alpha) \log Y }{ \alpha \log Y - \log h }. 
\end{split}\label{splitloglog}
\end{align}
 We then find that \eqref{loghleq} and \eqref{splitloglog} together give us that 
\begin{align}
\begin{split}
 \frac{ \log\left( \frac{Y^{\delta}}{h} \right) }{ \log\left( \frac{Y^{\alpha}}{h} \right) }
 & \leq 1 + \frac{\delta - \alpha}{\alpha - \gamma} \\ 
 & = \frac{\delta - \gamma}{\alpha - \gamma}. 
\end{split}\label{simplifyloglog}
\end{align}
 From \eqref{simplifyloglog} and the above bound for the given sum of reciprocals of primes, we find that 
\begin{equation}\label{miscsimplify}
 \sum_{\substack{ Y^{\alpha} < p \leq Y^{\delta} \\ p \equiv b \bmod {h} }} 
 \frac{1}{p} 
 \leq \frac{2}{\varphi(h)} 
 \left( \frac{1}{\log\left( \frac{Y^{\delta}}{h} \right)} + \log\left( \frac{\delta - \gamma}{\alpha - \gamma} \right) \right). 
\end{equation}
 Now, from \eqref{loghleq}, we find that 
\begin{equation}\label{logYpoweroverh} 
 \log\left( \frac{Y^{\delta}}{h} \right) \geq (\delta - \gamma) \log Y. 
\end{equation}
 A combined application of \eqref{Ygeq2}, \eqref{miscsimplify}, and 
 \eqref{logYpoweroverh} then gives us the desired result. 
 \end{proof}

\begin{lemma}
 Let $\alpha$ and $\delta$ be fixed, with $0 < \alpha < \delta$, and let $\beta, K > 0$ be fixed. Also, let $d$ be a fixed nonnegative integer, 
 and let $M > 0$ be fixed. For each real $Y \geq 2$, 
 let $P(T) = P = P_Y \in \mathbb{Z}[T]$ denote a nonzero polynomial of degree at most $d$ possibly depending on $Y$ 
 such that 
\begin{equation}\label{absPt} 
 |P(t)| \leq Y^{M} 
\end{equation}
 for $Y^{\alpha} < t \leq Y^{\delta}$. Then 
\begin{equation}\label{llmany} 
 \sum_{\substack{ Y^{\alpha} < p \leq Y^{\delta} \\ P(p) \neq 0 }} 
 \frac{ \mathcal{B}_{\beta, K}(|P(p)|) }{p}
 \ll_{\alpha, \delta, \beta, K, d, M} \, \mathcal{B}_{\beta, K}(\operatorname{cont}(P)), 
\end{equation}
 as $Y \to \infty$, uniformly in $P$. 
\end{lemma}

\begin{proof}
 We write $P = \operatorname{cont}(P) Q$, where $Q$ is primitive. Observe that $Q \in \mathbb{Z}[T]$ and that $\deg Q = \deg P 
 \leq d$, and observe that 
\begin{align}
\begin{split}
 |Q(u)| 
 & = \frac{|P(u)|}{\operatorname{cont}(P)} \\ 
 & \leq |P(u)| 
\end{split}\label{splituseabs}
\end{align} 
 for $u \in (Y^{\alpha}, Y^{\delta}]$.
 From \eqref{splituseabs} and the assumption in \eqref{absPt}, 
 we see that 
\begin{equation}\label{QuleqYM}
 |Q(u)| \leq Y^{M} 
\end{equation}
 for $u \in (Y^{\alpha}, Y^{\delta}]$.
 Thus $P$ and $Q$ satisfy the same required conditions specified above. 

 Since $\operatorname{cont}(P) > 0$, we find that $P(p) \neq 0 \Longleftrightarrow Q(p) \neq 0$ for a given prime $p$, 
 and thus the sum on the left of 
 \eqref{llmany} runs over the same primes
 as the corresponding sum obtained
 by replacing the condition $P(p) \neq 0$ with $Q(p) \neq 0$. 
 So, for a prime $p$ in either sum, the value 
 $|Q(p)|$ is a positive integer, with 
\begin{equation}\label{absPp} 
 |P(p)| = \operatorname{cont}(P) \, |Q(p)|. 
\end{equation}
 From \eqref{absPp} and Lemma \ref{lemcalBineq}, we find that 
\begin{align}
\begin{split}
 \mathcal{B}_{\beta, K}(|P(p)|)
 & = \mathcal{B}_{\beta, K}( \operatorname{cont}(P) \, |Q(p)| ) \\ 
 & \leq \mathcal{B}_{\beta, K}(\operatorname{cont}(P)) \, \mathcal{B}_{\beta, K}(|Q(p)|).
\end{split}\label{splitcalBcont}
\end{align} 
 From \eqref{splitcalBcont} (recalling a previous note on the index set of the sum in \eqref{llmany}), 
 we obtain that 
\begin{equation}\label{reducetoQ}
 \sum_{\substack{ Y^{\alpha} < p \leq Y^{\delta} \\ P(p) \neq 0 }} 
 \frac{ \mathcal{B}_{\beta, K}(|P(p)|) }{p} 
 \leq \mathcal{B}_{\beta, K}(\operatorname{cont}(P))
 \sum_{\substack{ Y^{\alpha} < p \leq Y^{\delta} \\ Q(p) \neq 0 }} 
 \frac{ \mathcal{B}_{\beta, K}(|Q(p)|) }{p}. 
\end{equation}
 Now, for a primitive polynomial $Q$, we have that $\operatorname{cont}(Q) = 1$ and that 
 $\mathcal{B}_{\beta, K}(\operatorname{cont}(Q)) = \mathcal{B}_{\beta, K}(1) = 1$. So, from \eqref{reducetoQ}, 
 we find that the desired result holds if and only if 
 $$ \sum_{\substack{ Y^{\alpha} < p \leq Y^{\delta} \\ Q(p) \neq 0 }} 
 \frac{ \mathcal{B}_{\beta, K}(|Q(p)|) }{p}
 \ll_{\alpha, \delta, \beta, K, d, M} 1 $$
 holds uniformly for every nonzero primitive polynomial $Q \in \mathbb{Z}[T]$ 
 of degree at most $d$ satisfying $|Q(u)| \leq Y^{M}$ for every real $u \in (Y^{\alpha}, Y^{\delta}]$.

 Expanding the product in \eqref{definecalB}, we find that 
\begin{equation}\label{eq:Bexpand}
 \mathcal{B}_{\beta,K}(n)
 =\sum_{\substack{m \mid n\\ m \ \text{squarefree}}}
 \frac{K^{\omega(m)}}{m^\beta} 
\end{equation}
 for $n \geq 1$ (noting that the sum in \eqref{eq:Bexpand} is not necessarily only over prime divisors).
 Now, fix a real number 
\begin{equation}\label{boundtheta} 
 \theta \in (0, \alpha), 
\end{equation}
 with the choice of $\theta$ depending only on $\alpha$. Also, consider the index set associated with the sum in \eqref{eq:Bexpand} 
 (consisting of positive and squarefree 
 integers $m$
 such that $m \mid n$), and 
 consider indices of this form such that $m > Y^{\theta}$. 
 Since $m > Y^{\theta}$ in this case, we find that 
 $ \frac{1}{m^{\beta}} \leq \left( Y^{\theta} \right)^{-\frac{\beta}{2}} \frac{1}{m^{\frac{\beta}{2}}}$, 
 so that 
\begin{align}
\begin{split}
 \sum_{\substack{m \mid n \\ m \ \text{squarefree} \\ m > Y^{\theta}}}
 \frac{K^{\omega(m)}}{m^{\beta}} 
 & \leq Y^{-\frac{\theta \beta}{2}} \sum_{\substack{m \mid n \\ m \ \text{squarefree} \\ 
 m > Y^{\theta} }} \frac{K^{\omega(m)}}{m^{\frac{\beta}{2}}} \\ 
 & \leq Y^{-\frac{\theta \beta}{2}} \sum_{\substack{m \mid n \\ m \ \text{squarefree}}} 
 \frac{K^{\omega(m)}}{m^{\frac{\beta}{2}}}.
\end{split}\label{splittopower}
\end{align}
 Now, observe that 
\begin{align*}
 \sum_{\substack{m \mid n \\ m \ \text{squarefree}}} 
 \frac{K^{\omega(m)}}{m^{\frac{\beta}{2}}}
 & = \prod_{\substack{ p \mid n \\ p \ \text{prime}}} \left( 1 + \frac{K}{p^{\frac{\beta}{2}}} \right) \\ 
 & \leq \prod_{\substack{p \mid n \\ p \ \text{prime}}} (1 + K) \\ 
 & = (1 + K)^{\omega(n)}. 
\end{align*}
 From this and \eqref{splittopower} together, we find that 
\begin{equation}\label{wanto1} 
 Y^{-\frac{\theta \beta}{2}} 
 \sum_{\substack{m \mid n \\ m \ \text{squarefree}}} 
 \frac{K^{\omega(m)}}{m^{\frac{\beta}{2}}}
 \leq Y^{-\frac{\theta \beta}{2}} (1 + K)^{\omega(n)}. 
\end{equation}
 
 Now, rewrite the right-hand side of \eqref{wanto1} so that 
\begin{equation}\label{rewriteerror}
 Y^{-\frac{\theta \beta}{2}} (1+K)^{\omega(n)}
 = Y^{-\frac{\theta \beta}{2} + \varepsilon(Y, n)}, 
\end{equation}
 where the error term $\varepsilon(Y, n)$ in \eqref{rewriteerror} is such that 
\begin{equation}\label{defineerror}
 \varepsilon(Y, n) = \frac{\omega(n) \log(1+K)}{\log Y}. 
\end{equation}
 We seek, at this point, to estimate the error term defined in 
 \eqref{defineerror}. 

 Suppose that $Y > \exp(1)$. We partition the distinct prime 
 divisors of $n$ according to whether such factors are $\leq \sqrt{\log Y}$ or $> \sqrt{\log Y}$. Of course, there are at most $\sqrt{\log Y}$ 
 prime divisors in this first part of the given partition (i.e., since there are at most $\sqrt{\log Y}$ positive integers $\leq \sqrt{\log Y}$). 
 To estimate the number of prime divisors in the second part, we now restrict our attention to positive integers $n \leq Y^{M}$. Now, 
 let $k$ denote the number of distinct prime divisors of $n$ exceeding 
 $\sqrt{\log Y}$. Their product divides $n$, so that 
\begin{equation}\label{rootlogYk} 
 \left( \sqrt{\log Y} \right)^{k} \leq n \leq Y^{M}, 
\end{equation}
 and, from \eqref{rootlogYk}, we deduce that $ k \leq \frac{2 M \log Y}{\log \log Y}$. 
 Combining the given bounds on
 the cardinalities of the sets of prime divisors of $n$ that we have formed, 
 we find that 
\begin{equation}\label{partitionomega} 
 \omega(n) \leq \sqrt{\log Y} + \frac{2 M \log Y}{\log\log Y}
\end{equation}
 for each real $Y > \exp(1)$ and for each positive integer $n \leq Y^{M}$. Since $K$ and $M$ are fixed, the definition of $\varepsilon(Y, n)$ in 
 \eqref{defineerror} and the bound in \eqref{partitionomega} together give us that 
\begin{align}
\begin{split}\label{boundsup}
 0 & \leq \sup_{\substack{n\in\mathbb Z\\1\leq n\leq Y^M}}
\varepsilon(Y,n) \\ 
 & \leq \log(1 + K) \left( \frac{1}{\sqrt{\log Y}} + \frac{2 M}{\log \log Y} \right). 
\end{split}
\end{align}
 From \eqref{boundsup}, we find that 
\begin{equation}\label{limbracketsup}
 \lim_{Y \to \infty} \left( \sup_{\substack{n\in\mathbb Z\\1\leq n\leq Y^M}}
\varepsilon(Y,n) \right) = 0. 
\end{equation}
 Since $\theta \beta > 0$, we have that $\varepsilon(Y, n) \leq \frac{\theta \beta}{4}$ for all sufficiently large $Y$ and for every positive 
 integer $n \leq Y^{M}$. Now, the relations among 
 \eqref{splittopower} and \eqref{wanto1} and \eqref{rewriteerror} together give us that 
\begin{align}
\begin{split}
 0 & \leq \sum_{\substack{ m \mid n \\ m \ \text{squarefree} \\ m > Y^{\theta} }} 
 \frac{K^{\omega(m)}}{m^{\beta}} \\ 
 & \leq Y^{-\frac{\theta \beta}{2} + \varepsilon(Y, n)}. 
\end{split}\label{wantYvanish}
\end{align}
 So, from \eqref{limbracketsup} and \eqref{wantYvanish}, we find that 
\begin{equation}\label{tailto0}
 \sum_{\substack{ m \mid n \\ m \ \text{squarefree} \\ m > Y^{\theta} }} 
 \frac{K^{\omega(m)}}{m^{\beta}} \to 0 
\end{equation}
 as $Y \to \infty$, uniformly for integers $n \in [1, Y^{M}]$. 
 So, from \eqref{eq:Bexpand} and \eqref{tailto0} together, we find that 
\begin{equation}\label{calBo1}
 \mathcal{B}_{\beta, K}(n)
 = \sum_{\substack{m \mid n \\ m \ \text{squarefree} \\ m \leq Y^{\theta}}} 
 \frac{K^{\omega(m)}}{m^{\beta}} + o(1) 
\end{equation}
 as $Y \to \infty$, uniformly for integers $n \in [1, Y^{M}]$. 

 Recall that if $Q(p) \neq 0$ for a prime $p \in (Y^{\alpha}, Y^{\delta}]$, then 
 $|Q(p)|$ is a positive integer. 
 So, from this and \eqref{QuleqYM} together, we find that: If $Q(p) \neq 0$, then 
 $1 \leq |Q(p)| \leq Y^{M}$ whenever $p \in (Y^\alpha, Y^{\delta}]$. 
 So, from this and the estimate in \eqref{calBo1}, we find that 
\begin{equation}\label{multby1overp} 
 \mathcal{B}_{\beta, K}( |Q(p)| ) = \sum_{\substack{ m \mid |Q(p)| 
 \\ m \ \text{squarefree} \\ m \leq Y^{\theta} }} 
 \frac{K^{\omega(m)}}{m^{\beta}} + o(1) 
\end{equation}
 as $Y \to \infty$, uniformly over all nonzero primitive polynomials $Q \in \mathbb{Z}[T]$ 
 of degree at most $d$ satisfying $|Q(u)| \leq Y^{M}$ for each real value $u \in (Y^{\alpha}, Y^{\delta}]$,
 and all primes $p \in (Y^\alpha, Y^\delta]$
 for which $Q(p) \neq 0$. 

 Now, by taking $h = b = 1$ and $\gamma = \theta \in (0, \alpha)$ in 
 Lemma \ref{lemreciprocal}, this same lemma gives us that 
\begin{equation}\label{controlerror}
 \sum_{Y^{\alpha} < p \leq Y^{\delta}} \frac{1}{p} \ll_{\alpha, \delta} 1 
\end{equation}
 for all real $Y \geq 2$. Multiplying both sides of \eqref{multby1overp} by $\frac{1}{p}$ (for a prime $p$ such that $Y^{\alpha} < p 
 \leq Y^{\delta}$ and such that $Q(p) \neq 0$), 
 and then applying the summation operator 
 $\sum_{\substack{Y^{\alpha} < p \leq Y^{\delta} \\ Q(p) \neq 0}} \cdot$ to both sides
 of the resultant identity, 
 and by then exploiting the relation in 
 \eqref{controlerror} to control the resulting error term, we find that 
\begin{equation}\label{introducedouble}
 \sum_{\substack{ Y^{\alpha} < p \leq Y^{\delta} \\ Q(p) \neq 0 }} 
 \frac{\mathcal{B}_{\beta, K}(|Q(p)|)}{p} 
 = \sum_{\substack{ Y^{\alpha} < p \leq Y^{\delta} \\ Q(p) \neq 0 }}
 \frac{1}{p} \sum_{\substack{ m \mid |Q(p)| \\ m \ \text{squarefree} \\ m \leq Y^{\theta} }} 
 \frac{K^{\omega(m)}}{m^{\beta}} + o(1) 
\end{equation}
 as $Y \to \infty$, uniformly over all nonzero primitive polynomials $Q \in \mathbb{Z}[T]$ of degree at most $d$ such that $|Q(u)| \leq 
 Y^{M}$ for each real value $u \in (Y^{\alpha}, Y^{\delta}]$. Interchanging the order of summation in the double sum on the right of 
 \eqref{introducedouble}, we find that 
\begin{equation}\label{interchangeddouble} 
 \sum_{\substack{ Y^{\alpha} < p \leq Y^{\delta} \\ Q(p) \neq 0 }} 
 \frac{ \mathcal{B}_{\beta, K}(|Q(p)|) }{p}
 = \sum_{\substack{ 1 \leq m \leq Y^{\theta} \\ m \ \text{squarefree} }}
 \frac{K^{\omega(m)}}{m^{\beta}} 
 \sum_{\substack{ Y^{\alpha} < p \leq Y^{\delta} \\ 
 Q(p) \neq 0 \\ m \mid |Q(p)| }}
 \frac{1}{p} + o(1) 
\end{equation}
 as $Y \to \infty$, with the same uniformity in $Q$ as before. Now, recall the given bounds on $\theta$ indicated in \eqref{boundtheta}. 
 We thus find that each pair $(m, p)$ arising in the double sum in \eqref{interchangeddouble} is such that 
\begin{equation}\label{mlessp} 
 1 \leq m \leq Y^{\theta} < Y^{\alpha} < p. 
\end{equation}
 From \eqref{mlessp}, since $p$ is prime, we deduce that $\gcd(p, m) = 1$. Now, for an arbitrary positive integer $m$, define
 $$ \rho_{Q}(m) = \big| \{ b \in \mathbb{Z} : 1 \leq b \leq m, 
 \ \gcd(b, m) =1, \ Q(b) \equiv 0 \bmod {m} \} \big|. $$
 For each positive integer $m \leq Y^{\theta}$, by grouping the primes appearing 
 in the inner sum on the right of \eqref{interchangeddouble} 
 by their residue classes modulo $m$, 
 and by applying Lemma \ref{lemreciprocal}
 with $h = m$ and $\gamma = \theta$, we find that 
\begin{align}
\begin{split}
 \sum_{\substack{ Y^{\alpha} < p \leq Y^{\delta} \\ Q(p) \neq 0 \\ m \mid |Q(p)| }}
 \frac{1}{p} 
 & \leq \sum_{\substack{ 1 \leq b \leq m \\ \gcd(b, m) = 1 \\ Q(b) \equiv 0 \bmod {m} }}
 \sum_{\substack{ Y^{\alpha} < p \leq Y^{\delta} \\ p \equiv b \bmod {m} }} 
 \frac{1}{p} \\ 
 & \ll_{\alpha, \delta, \theta} \frac{\rho_{Q}(m)}{\varphi(m)}. 
\end{split}\label{rhooverphi}
\end{align}
 Since $Q$ is primitive, its reduction modulo any prime $q$ is a nonzero polynomial of degree at most $d$
 over the field $\mathbb{Z} / q \mathbb{Z}$, and thus it has at most $d$ roots, with 
\begin{equation}\label{rholeqd} 
 \rho_{Q}(q) \leq d. 
\end{equation} 
 By the Chinese remainder theorem, we have that 
\begin{equation}\label{fromCRT} 
 \rho_{Q}(m) = \prod_{\substack{q \mid m \\ q \ \text{prime}}} \rho_{Q}(q)
\end{equation}
 for each positive squarefree integer $m$, with $\rho_{Q}(1) = 1$ corresponding to the empty product. 
 So, from \eqref{rholeqd} and \eqref{fromCRT} together, by fixing $D = \max\{ 1, d \}$, 
 we find that 
\begin{equation}\label{rhoDpower}
 \rho_{Q}(m) \leq D^{\omega(m)} 
\end{equation}
 for each positive and squarefree integer $m$. 

 So, a combined application of \eqref{interchangeddouble}, \eqref{rhooverphi}, and \eqref{rhoDpower} then gives us that 
\begin{equation}\label{introduceC}
 \sum_{\substack{ Y^{\alpha} < p \leq Y^{\delta} \\ Q(p) \neq 0 }} 
 \frac{ \mathcal{B}_{\beta, K}(|Q(p)|) }{p}
 \leq C_{\alpha, \delta, \theta} \sum_{\substack{1 \leq m \leq Y^{\theta} \\ m \ \text{squarefree} }} 
 \frac{ (KD)^{\omega(m)} }{m^{\beta} \varphi(m)} + o(1) 
\end{equation}
 as $Y \to \infty$, with the same uniformity in $Q$, where $C_{\alpha, \delta, \theta} > 0$ 
 depends only on $\alpha$, $\delta$, and $\theta$.

 By Euler's product formula in \eqref{Eulerprod}, we write
\begin{equation}\label{KDprod}
 \frac{ (KD)^{\omega(m)} }{ m^{\beta} \varphi(m) } 
 = \prod_{p \mid m} \frac{KD}{p^{\beta} (p-1)} 
\end{equation}
 for squarefree $m$. Now, define $R_{Y} = \prod_{p \leq Y^{\theta}} p$. Using \eqref{KDprod}, we find that 
\begin{align}
\begin{split}
 \prod_{p \leq Y^{\theta}} \left( 1 + \frac{KD}{p^{\beta}(p-1)} \right)
 & = \sum_{m \mid R_{Y}} \prod_{p \mid m} \frac{KD}{p^{\beta}(p-1)} \\ 
 & = \sum_{m \mid R_{Y}} \frac{(KD)^{\omega(m)}}{m^{\beta} \varphi(m)}. 
\end{split}\label{prodsumprod}
\end{align}
 If $m$ is squarefree and $1 \leq m \leq Y^{\theta}$, then each prime divisor $p$ of $m$ is such that 
 $p \leq m \leq Y^{\theta}$. 
 So, every prime divisor of $m$ occurs in $R_{Y}$.
 Since $m$ is squarefree, this implies that $m \mid R_{Y}$, and we thus find, using \eqref{prodsumprod}, that 
\begin{align}
\begin{split}
 \sum_{\substack{1 \leq m \leq Y^{\theta} \\ m \ \text{squarefree} }}
 \frac{ (KD)^{\omega(m)} }{m^{\beta} \varphi(m)}
 & \leq \sum_{m \mid R_{Y}} \frac{ (KD)^{\omega(m)} }{m^{\beta} \varphi(m)} \\ 
 & = \prod_{p \leq Y^{\theta}} \left( 1 + \frac{KD}{p^{\beta}(p-1)} \right). 
\end{split}\label{sumtoprod}
\end{align} 
 Using the inequality $1 + u \leq \exp(u)$ for $u \geq 0$, and by appealing to the inequality $p-1 \geq \frac{p}{2}$
 holding for each prime $p$, we find that 
\begin{align}
\begin{split}
 \prod_{p \leq Y^{\theta}} \left( 1 + \frac{KD}{p^{\beta} (p-1)} \right)
 & \leq \exp\left( KD \sum_{p \leq Y^{\theta}} \frac{1}{p^{\beta} (p-1)} \right) \\
 & \leq \exp\left( 2 K D \sum_{n=2}^{\infty} \frac{1}{n^{1 + \beta}} \right), 
\end{split}\label{prodtoexp}
\end{align} 
 and we see that the right-hand side of the latter inequality above is finite since $\beta > 0$, 
 and that this right-hand side depends only on $\beta$, $K$, and $d$, 
 recalling that $D = \max\{ 1, d \}$. 

 So, from a combined application of \eqref{introduceC}, \eqref{sumtoprod}, 
 and \eqref{prodtoexp}, for fixed $\theta \in (0, \alpha)$, we find that 
\begin{equation}\label{equivalentestimate} 
 \sum_{\substack{ Y^{\alpha} < p \leq Y^{\delta} \\ Q(p) \neq 0 }} 
 \frac{\mathcal{B}_{\beta, K}(|Q(p)|)}{p} \ll_{\alpha, \delta, \beta, K, d, M} 1 
\end{equation}
 as $Y \to \infty$, uniformly over all nonzero primitive polynomials $Q \in \mathbb{Z}[T]$
 of degree at most $d$ satisfying $|Q(u)| \leq Y^{M}$ for every real $u \in (Y^{\alpha}, Y^{\delta}]$. 
 Applying the estimate in \eqref{equivalentestimate}
 to $Q = \frac{P}{\operatorname{cont}(P)}$ in \eqref{reducetoQ}, we find that 
 $$ \sum_{\substack{Y^{\alpha} < p \leq Y^{\delta} \\ P(p) \neq 0}} 
 \frac{\mathcal{B}_{\beta, K}(|P(p)|)}{p} \ll_{\alpha, \delta, \beta, K, d, M} 
 \, \mathcal{B}_{\beta, K}(\operatorname{cont}(P)) $$
 as $Y \to \infty$, uniformly over all nonzero polynomials $P \in \mathbb{Z}[T]$
 satisfying the conditions of the lemma, thus completing our proof. 
\end{proof}

\begin{remark}
The constant implicit in \eqref{llmany} can be made explicit. With the parameters as in the preceding lemma, set
$$ D = \max\{1, d\}, \qquad
\eta=\frac{\alpha\beta}{4M},
\qquad
R=\max\left\{
2,\left(\frac{K}{\eta\log 2}\right)^{\frac{2}{\beta}}
\right\},
$$ and define the finite product $$ A = \prod_{p\leq R} \max\left\{
1,\frac{1+Kp^{-\frac{\beta}{2}}}{p^\eta}
\right\}.
$$
One may take
$$ 
 C = 
\left(
\frac{2}{ \left(\delta - \frac{\alpha}{2}\right)\log 2}
+
2\log\left(\frac{2\delta-\alpha}{\alpha}\right)
\right) 
\left(
\exp\left(\frac{2KD}{\beta}\right)+A
\right).
 $$
 For this choice of $C$, it can be shown (omitting details) that the inequality
\[
\sum_{\substack{Y^\alpha<p\leq Y^\delta\\P(p)\neq0}}
\frac{\mathcal B_{\beta,K}(|P(p)|)}p
\leq
C\,\mathcal B_{\beta,K}(\operatorname{cont}(P))
\]
holds for every real $Y\geq2$ and every polynomial $P$
satisfying the hypotheses of the preceding lemma.
\end{remark}

\begin{lemma}
 Let $\beta > 0$ and $K > 0$ be fixed. Then there exists a constant $C_{\beta, K} > 0$, depending only on $\beta$ and $K$, such that 
 $$ \sum_{\substack{1 \leq a \leq Z \\ 1 \leq b \leq Z}} \frac{ \mathcal{B}_{\beta, K}(\gcd(a, b)) }{ab} 
 \leq C_{\beta, K} (\log Z)^2 $$
 for each real value $Z \geq 2$, 
 where $a$ and $b$ are positive integers. Equivalently, the relation 
 $$ \sum_{\substack{1 \leq a \leq Z \\ 1 \leq b \leq Z}} 
 \frac{ \mathcal{B}_{\beta, K}(\gcd(a, b)) }{ab} 
 \ll_{\beta, K} (\log Z)^2 $$
 holds. 
 \end{lemma}

\begin{proof}
 Expanding the right-hand side of \eqref{definecalB}, we find that 
 $$ \mathcal{B}_{\beta, K}(\gcd(a, b))
 = \sum_{\substack{ m \mid a \\ m \mid b \\ m \ \text{squarefree} }} 
 \frac{K^{\omega(m)}}{m^{\beta}}. $$
 Using an interchange of finite summation operators, we find that 
\begin{align}
\begin{split}
 \sum_{\substack{ 1 \leq a \leq Z \\ 1 \leq b \leq Z }}
 \frac{\mathcal{B}_{\beta, K}(\gcd(a, b))}{ab}
 & = \sum_{\substack{1 \leq m \leq Z \\ m \ \text{squarefree}}}
 \frac{K^{\omega(m)}}{m^{\beta}} \left( \sum_{\substack{1 \leq a \leq Z \\ m \mid a}} \frac{1}{a} \right)^{2} \\
 & = \sum_{\substack{1 \leq m \leq Z \\ m \ \text{squarefree}}} 
 \frac{K^{\omega(m)}}{m^{\beta + 2}} \left( \sum_{1 \leq k \leq \frac{Z}{m}} \frac{1}{k} \right)^{2}.
\end{split}\label{introduceharmonic}
\end{align}
 A standard estimate for harmonic numbers gives us that 
 $\sum_{1 \leq k \leq \frac{Z}{m}} \frac{1}{k} \leq 1 + \log\left( \frac{Z}{m} \right) 
 \leq 1 + \log Z$. This and \eqref{introduceharmonic} together give us that 
 $$ \sum_{\substack{1 \leq a \leq Z \\ 1 \leq b \leq Z}}
 \frac{\mathcal{B}_{\beta, K}(\gcd(a, b))}{ab}
 \leq (1 + \log Z)^{2} \sum_{\substack{1 \leq m \leq Z \\ m \ \text{squarefree}}} 
 \frac{K^{\omega(m)}}{m^{\beta + 2}}. $$
 Direct analogues of the derivation of the ``sum-to-product'' 
 comparison in \eqref{sumtoprod} and subsequent arguments give us that 
\begin{align*}
 \sum_{\substack{1 \leq m \leq Z \\ m \ \text{squarefree}}} \frac{K^{\omega(m)}}{m^{\beta + 2}}
 & \leq \prod_{p \leq Z} \left( 1 + \frac{K}{p^{\beta +2}} \right) \\
 & \leq \exp\left( K \sum_{p \leq Z} \frac{1}{p^{\beta + 2}} \right) \\
 & \leq \exp\left( K \sum_{n=2}^{\infty} \frac{1}{n^{\beta + 2}} \right) \\ 
 & \leq \exp\left( \frac{K}{\beta + 1} \right), 
\end{align*}
 where the last inequality follows 
 from $\sum_{n=2}^{\infty} \frac{1}{n^{\beta + 2}} \leq \int_{1}^{\infty} \frac{dt}{t^{\beta + 2}} 
 = \frac{1}{\beta + 1}$. 
 Finally, since 
 $Z \geq 2 \Longrightarrow 1 + \log Z \leq \big( 1 + \frac{1}{\log 2} \big) \log Z$, 
 we find that a permissible choice for the constant $C_{\beta, K}$
 is $C_{\beta, K} = \big( 1 + \frac{1}{\log 2} \big)^{2} \exp\big( \frac{K}{\beta + 1} \big)$. 
\end{proof}

 If $a$ and $b$ are positive integers and $p$ and $q$ are primes such that $p \nmid a$ and $q 
 \nmid b$, Lemma \ref{lemmawap} gives us that 
\begin{equation}\label{wlinear} 
 w(ap) = w(bq)
 \Longleftrightarrow w(a) p - w(b) q = b - a. 
\end{equation}
 Because of this, we need a uniform upper bound for the number of ordered
 prime pairs satisfying a linear equation as in \eqref{wlinear}. 

 The proof of Lemma~\ref{decimallemma} below uses Selberg's upper bound sieve method; see \cite[Section~7.2]{CojocaruMurty2005}
 for a general treatment. In that proof, we bound $S$ by a sum of
squares of weighted divisor sums and choose the weights
$\lambda_h$ so that the resulting quadratic form equals
$\frac{1}{\mathcal{G}(z)}$, leading to \eqref{Sboundx003}.
We give the argument in detail for the two linear polynomials
in \eqref{linearT}, including the estimates needed to obtain
a bound uniform in $a$, $b$, $A$, and $B$.

\begin{lemma}\label{decimallemma}
 There exist absolute constants $C > 0$ and 
\begin{equation}\label{x0geq8}
 x_0 \geq 8
\end{equation}
 such that the following holds. Let 
\begin{equation}\label{xgeqx0} 
 x \geq x_0 
\end{equation} 
 be real, and let $a$, $b$, $A$, 
 and $B$ be positive integers such that 
\begin{align}
 a & \neq b, \nonumber \\ 
 1 & \leq a, b \leq x^{0.054}, \label{abinterval} \\ 
 \frac{a}{2} & \leq A \leq x^{0.055}, \text{and} \label{Ainterval} \\ 
 \frac{b}{2} & \leq B \leq x^{0.055}. \label{bover2leqB}
\end{align}
 Then the number of ordered pairs $(p, q)$ of primes satisfying
\begin{align}
 A p - B q & = b - a, \nonumber \\ 
 \frac{x}{2a} & < p \leq \frac{x}{a}, \text{and} \label{pinterval} \\
 \frac{x}{2b} & < q \leq \frac{x}{b} \label{qinterval}
\end{align}
 is $0$ if $\gcd(A, B) \nmid (b-a)$ and is at most 
\begin{equation}\label{miscbound}
 C \left( \frac{x \, \gcd(A, B) }{a b (\log x)^2} \mathcal{B}_{1,2}(AB |a-b|) + 
 x^{0.03} \right) 
\end{equation}
 if $\gcd(A, B)$ divides $b - a$. 
 The same constants $C$ and $x_0$ work for each choice of $a$, $b$, $A$, and $B$
 satisfying the stated conditions (and where the integers $a$, $b$, $A$, and $B$ may depend on $x$). 
\end{lemma}

\begin{proof}
 Since $\gcd(A, B)$ divides $Au - Bv$ for all integers $u$, $v$, the equation $A u - B v = b - a$ has no solutions when $\gcd(A, B) \nmid (b 
 - a)$. So, we henceforth assume that $\gcd(A, B)$ divides $b-a$. 
 
 By B\'{e}zout's identity (recalling Section \ref{secprelim}) together with 
 the assumption that $\gcd(A, B) \mid (b-a)$, 
 we are permitted to fix integers $u_0$ and $v_0$ 
 such that 
\begin{equation}\label{fromBezout}
 A u_0 - B v_0 = b - a. 
\end{equation}
 Since $\frac{A}{\gcd(A, B)}$ and $\frac{B}{\gcd(A, B)}$ are coprime, 
 all integer solutions of $A u - B v = b - a$ are parametrized uniquely by 
\begin{equation}\label{rewritewithp}
 u = u_0 + \frac{B}{\gcd(A, B)} k 
\end{equation} 
 and 
\begin{equation}\label{rewritewithq} 
 v = v_0 + \frac{A}{\gcd(A, B)} k 
\end{equation}
 for $k \in \mathbb{Z}$.
 Each pair $(p, q)$ of primes counted in the lemma corresponds to one 
 of these integer solutions with $u = p$ and $v = q$ 
 and hence to an integer $k$ such that 
\begin{equation}\label{Bovergcd}
 \frac{x}{2a} < u_0 + \frac{B}{\gcd(A, B)} k \leq \frac{x}{a} 
\end{equation}
 and $$ \frac{x}{2b} < v_0 + \frac{A}{\gcd(A, B)} k \leq \frac{x}{b}. $$
 The relations in \eqref{Bovergcd} give us that 
\begin{equation}\label{intervalbounds}
 \frac{\gcd(A, B)}{B}\left( \frac{x}{2a} - u_0 \right) 
 \leq k \leq \frac{\gcd(A, B)}{B}\left( \frac{x}{a} - u_{0}	 \right), 
\end{equation}
 so the corresponding integers $k$ lie in an interval of length
 $$ \frac{x \gcd(A, B)}{2a B} \leq \frac{x \, \gcd(A, B)}{ab}, $$
 where we have used the relation $B \geq \frac{b}{2}$ appearing in 
 \eqref{bover2leqB}. 

 As above, we assume that
 $x_{0} \geq 8$. From the relations in \eqref{pinterval} and \eqref{qinterval}, along with the bounds $a, b \leq x^{0.054}$ 
 and $A, B \leq x^{0.055}$ given in \eqref{abinterval}, \eqref{Ainterval}, and \eqref{bover2leqB}, we
 find that 
\begin{align*}
 \min\{ p, q \}
 & > \frac{x}{2 \max\{ a, b \}} \\ 
 & \geq \frac{1}{2} x^{0.946}, 
\end{align*}
 for each prime pair $(p, q)$ satisfying the conditions of the lemma. 
 In turn, this gives us that 
\begin{align}
\begin{split}
 \min\{ p, q \} & > x^{0.055} \\
 & \geq \max\{ A, B \}. 
\end{split}\label{splitsplit}
\end{align}
 Recalling that we are letting $u = p$ and $v = q$, we obtain from 
 \eqref{rewritewithp} and \eqref{rewritewithq} that 
\begin{equation}\label{preplacesu}
 p = u_0 + \frac{B}{\gcd(A, B)} k 
\end{equation} 
 and 
\begin{equation}\label{qreplacesv}
 q = v_0 + \frac{A}{\gcd(A, B)} k. 
\end{equation}
 Now, observe that $\gcd\big( u_0, \frac{B}{\gcd(A, B)} \big)$ divides both sides of \eqref{preplacesu} and that $\gcd\big( v_0, 
 \frac{A}{\gcd(A, B)} \big)$ divides both sides of \eqref{qreplacesv}. Since $p$ and $q$ are prime, and since 
 $$ 1 \leq \gcd\left( u_{0}, \frac{B}{\gcd(A, B)} \right) \leq B < p$$ and 
 $$ 1 \leq \gcd\left( v_{0}, \frac{A}{\gcd(A, B)} \right) \leq A < q, $$
 the divisibility relations just established give us that a prime pair
 satisfying the lemma's conditions can exist only if 
\begin{equation}\label{gcdgcd1}
 \gcd\left( u_{0}, \frac{B}{\gcd(A, B)} \right) = \gcd\left( v_{0}, \frac{A}{\gcd(A, B)} \right) = 1.
\end{equation}
 So, if either of 
 the linear polynomials in $\mathbb{Z}[T]$ among 
\begin{equation}\label{linearT} 
 u_0 + \frac{B}{\gcd(A, B)} T \quad \text{and} \quad v_{0} + \frac{A}{\gcd(A, B)} T 
\end{equation}
 is \emph{not} primitive, then (recalling the preceding result on the equalities in \eqref{gcdgcd1}) 
 the number of prime pairs satisfying the lemma's conditions is $0$.
 So, it remains to prove the bound involving \eqref{miscbound} when each of the polynomials
 indicated in \eqref{linearT} is primitive. 
 So, we henceforth assume that each of the polynomials in \eqref{linearT} is primitive. 

 For each positive integer $m$, let $\rho(m)$ denote the number of residue
 classes $k \bmod m$ satisfying
\begin{equation}\label{countrho}
 \left( u_{0} + \frac{B}{\gcd(A, B)} k \right) \left( v_{0} + \frac{A}{\gcd(A, B)} k \right) 
 \equiv 0 \pmod {m}.
\end{equation}
 From our assumption that each of the polynomials in \eqref{linearT} is primitive, we find that their reductions modulo any prime $p$ are 
 nonzero polynomials of degree at most $1$.
 Now, observe that the defining congruence in \eqref{countrho}, for $m = p$, 
 holds precisely when 
\begin{equation}\label{orcongruence1}
 \left( u_{0} + \frac{B}{\gcd(A, B)} k \right) \equiv 0 \bmod {p} 
\end{equation}
 or
\begin{equation}\label{orcongruence2}
 \left( v_{0} + \frac{A}{\gcd(A, B)} k \right) \equiv 0 \bmod{p}. 
\end{equation} 
 Since each of the reduced versions, modulo $p$, 
 of the polynomials among \eqref{linearT} is nonzero and has degree at most $1$
 over $\mathbb{Z} / p\mathbb{Z}$, 
 each of the two congruences among \eqref{orcongruence1} and 
 \eqref{orcongruence2} is satisfied by at most one residue class modulo $p$, 
 so that 
\begin{equation}\label{rhopleq2}
 \rho(p) \leq 2
\end{equation}
 for each prime $p$. 

 If $\rho(2) = 2$, then, for each integer $k$, 
 at least one of 
 $u_0 + \frac{B}{\gcd(A, B)} k$ and $v_0 + \frac{A}{\gcd(A, B)} k$ is even. 
 By way of contradiction, suppose that an ordered pair $(p, q)$ of primes
 satisfies the lemma's conditions. Now, consider the integer parameter $k$ associated with $(p, q)$
 according to \eqref{preplacesu} and \eqref{qreplacesv}. 
 For the integer $k$ associated with $(p, q)$, 
 the two linear expressions equal $p$ and $q$; 
 since $\rho(2) = 2$, at least one of these integers is even, so $2 \mid pq$. 
 Consequently, we find that $p$ or $q$ is even. 
 So, we have that $p = 2$ or $q = 2$, and thus $\min\{ p, q\} = 2$. 
 On the other hand, the condition $p > \frac{x}{2a}$ in 
 \eqref{pinterval}, the condition $q > \frac{x}{2b}$ 
 in \eqref{qinterval}, 
 the condition $a, b \leq x^{0.054}$ in \eqref{abinterval}, 
 and the condition $x \geq x_0 \geq 8$ given by 
 \eqref{x0geq8} and \eqref{xgeqx0} together give us that 
\begin{align*}
 \min\{ p, q \} 
 & > \frac{x}{2 \max\{ a, b \}} \\
 & \geq \frac{1}{2} x^{0.946} \\ 
 & \geq \frac{1}{2} 8^{0.946} \\
 & > 2.
\end{align*}
 This contradicts that $\min\{ p, q \} = 2$. So, when $\rho(2) = 2$, no ordered pair of primes satisfies the lemma's conditions 
 for $x \geq x_0 \geq 8$ (as required by \eqref{x0geq8} and \eqref{xgeqx0}). 

 So, since $\rho(2) \leq 2$, we may henceforth assume that $\rho(2) \leq 1$, since the case $\rho(2) =2$ gives zero prime pairs 
 satisfying the lemma's
 conditions. Now, recall that \eqref{rhopleq2} holds for each prime $p$. 
 So, since $\rho(p) \leq 2 < p$ for each odd prime $p$, 
 this and a previous argument together give us that 
\begin{equation}\label{rhostrictp}
 0\leq \rho(p) < p 
\end{equation}
 for each prime $p$. 

 For the time being, suppose that $p$ is a prime such that $p \nmid 2 A B |a - b|$. Thus $p\nmid A$ and $p \nmid B$. If $p$ were to 
 divide $\frac{B}{\gcd(A, B)}$, then $p$ would divide $\gcd(A, B) \frac{B}{\gcd(A, B)}$, contradicting that $p \nmid B$. So, again with 
 the assumption that $p \nmid 2 A B |a - b|$, 
 each of the integers among $\frac{B}{\gcd(A, B)}$ and $\frac{A}{\gcd(A, B)}$ 
 is invertible modulo $p$. So, each of the congruences among 
 \eqref{orcongruence1} and \eqref{orcongruence2} has exactly one solution modulo $p$. 
 The two congruences among \eqref{orcongruence1} and \eqref{orcongruence2} 
 have distinct solutions modulo $p$, since a common solution $k \bmod {p}$ would, 
 recalling the consequence of 
 B\'{e}zout's identity in \eqref{fromBezout}, 
 imply that 
\begin{align*}
 b - a 
 & = A u_{0} - B v_{0} \\
 & = A \left( u_{0} + \frac{B}{\gcd(A, B)} k \right) - B \left( v_{0} + \frac{A}{\gcd(A, B)} k \right) \\ 
 & \equiv 0 \pmod {p}, 
\end{align*}
 but this would contradict that $p \nmid |a - b|$. 

 By Euclid's lemma, the congruence defining $\rho(p)$, i.e., the $m = p$ case of \eqref{countrho}, 
 holds precisely when at least one of the congruences among
 \eqref{orcongruence1} and \eqref{orcongruence2} holds. 
 As demonstrated above, if $p$ is a prime such that $p \nmid 2 A B |a - b|$, 
 then each of the congruences among 
 \eqref{orcongruence1} and \eqref{orcongruence2} has exactly one solution modulo $p$, 
 and these solutions are distinct. So, for a prime $p$, we obtain the implication such that 
 $$ p \nmid 2 A B |a-b| \Longrightarrow \rho(p) = 2. $$

 By the Chinese remainder theorem, for each positive and squarefree integer $m$, we have that 
\begin{equation}\label{furtherCRT} 
 \rho(m) = \prod_{p \mid m} \rho(p), 
\end{equation}
 with the understanding that $\rho(1) = 1$ 
 and with the empty product being interpreted as $1$. 

 Now, let 
\begin{equation}\label{zx0005}
 z = x^{0.005}. 
\end{equation}
 So, each ordered 2-tuple $(p, q)$ of primes satisfying the lemma's conditions is such that 
\begin{equation}\label{introducedecimal}
 \min\{ p, q \} > x^{0.055} > z,
\end{equation}
 where the initial inequality in \eqref{introducedecimal}
 was established in \eqref{splitsplit}, 
 and where the latter inequality in 
 \eqref{introducedecimal} follows from the assumptions 
 such that $x \geq x_0 \geq 8$. 

 Recalling \eqref{intervalbounds}, we let $S$ denote the number of integers $k$ in the interval
\begin{equation}\label{intervalforS} 
 \left[ \frac{\gcd(A, B)}{B}\left( \frac{x}{2a} - u_{0} \right), \frac{\gcd(A, B)}{B}\left( \frac{x}{a} - u_{0} \right) \right] 
\end{equation}
 such that the integer $$ \left( u_{0} + \frac{B}{\gcd(A, B)} k \right) \left( v_{0} + \frac{A}{\gcd(A, B)}k \right) $$
 has no prime divisor at most $z$. So, again recalling \eqref{intervalbounds}, 
 according to the given definition of $S$, the number of ordered prime pairs
 satisfying the lemma's conditions is at most $S$. 
 Observe that the interval in \eqref{intervalforS}
 is of length $\frac{x \gcd(A, B)}{2 a B}$. Consequently, the number of integers in any fixed
 residue class modulo $m$ lying 
 in the interval in \eqref{intervalforS} differs from
 $\frac{x \gcd(A, B)}{2 a B m}$ by at most $1$. 
 So, the number of integers $k$ in the interval in \eqref{intervalforS}
 satisfying the congruence in \eqref{countrho}
 differs from
 $\frac{x \gcd(A, B) \rho(m)}{2a B m}$ by at most $\rho(m)$. 

 Define 
\begin{equation}\label{definecalG} 
 \mathcal{G}(z) = \sum_{\substack{1 \leq h \leq z \\ h \ \text{squarefree}}} \prod_{p \mid h} \frac{\rho(p)}{p- \rho(p)}, 
\end{equation}
 where the term corresponding to $h = 1$ is interpreted as $1$. From the inequalities shown in \eqref{rhostrictp}, we find that each 
 term in the sum in \eqref{definecalG} is nonnegative, 
 and, moreover, since the term corresponding to $h = 1$ equals $1$, 
 we find that $ \mathcal{G}(z) \geq 1$. 
 For a positive integer $h$, define $\lambda_h$
 so that $\lambda_1 = 1$ and so that $\lambda_h = 0$ whenever $h > z$ or $h$ is not squarefree, 
 and, for remaining values of $h$, let $\lambda_h$ (for the time being) 
 be understood as a fixed real number (to be made explicit later). 

 For every integer $k$ counted by $S$, we have 
\begin{equation}\label{contributes1}
 \sum_{\substack{1 \leq h \leq z \\ h \mid \big( u_0 + \frac{B}{\gcd(A, B)} k \big) 
 \big( v_0 + \frac{A}{\gcd(A, B)} k\big) }} \lambda_h = \lambda_1 = 1, 
\end{equation}
 since any index $h > 1$ in the sum would have a prime divisor $p \leq h \leq z$, 
 and the relations
 among $p \mid h$ and $h \mid \big( u_0 + \frac{B}{\gcd(A, B)} k \big) \big( v_{0} + \frac{A}{\gcd(A, B)} k \big)$
 would imply 
 $p \mid \big( u_0 + \frac{B}{\gcd(A, B)} k \big) \big( v_{0} + \frac{A}{\gcd(A, B)} k \big)$, 
 contradicting the requirement in the definition of $S$ that no prime at most $z$ can divide this product. 
 Since each integer $k$ counted by $S$ contributes $1$, according to \eqref{contributes1}, to the following sum 
 and since the remaining terms are squares of real numbers, we find that 
\begin{equation}\label{bigsquaresum} 
 S \leq \sum_{\substack{ k \in \mathbb{Z} \\ \frac{\gcd(A, B)}{B} \big( \frac{x}{2a} - u_0 \big)
 \leq k \\ k \leq \frac{\gcd(A, B)}{B} \big( \frac{x}{a} - u_0 \big) }} 
 \left( \sum_{\substack{1 \leq h \leq z \\ h \mid \big( 
 u_0 + \frac{B}{\gcd(A, B)} k \big) \big( v_0 + 
 \frac{A}{\gcd(A, B)} k \big) }} \lambda_h \right)^2.
\end{equation}
 Expanding the square in \eqref{bigsquaresum} (as a double sum) and interchanging the order of finite summation, we find that 
\begin{equation}\label{expandsquare}
S \leq
\sum_{1 \leq h,j \leq z} \lambda_h \lambda_j
\sum_{\substack{
k \in \mathbb{Z} \\[2pt]
\frac{\gcd(A,B)}{B}\left(\frac{x}{2a}-u_0\right) \leq k \\[2pt]
k \leq \frac{\gcd(A,B)}{B}\left(\frac{x}{a}-u_0\right) \\[2pt]
\operatorname{lcm}(h,j) \mid
\left(u_0+\frac{B}{\gcd(A,B)}k\right)
\left(v_0+\frac{A}{\gcd(A,B)}k\right)
}}
1,
\end{equation}
 where $h$ and $j$ range over positive integers. 
 Fix positive integers $h, j \leq z$.
 The divisibility condition associated with the inner sum 
 in \eqref{expandsquare} is equivalent to 
\begin{equation}\label{modlcm}
 \left( u_{0} + \frac{B}{\gcd(A, B)} k \right) \left( v_0 + \frac{A}{\gcd(A, B)} k \right) \equiv 
 0 \pmod { \operatorname{lcm}(h, j)}. 
\end{equation}
 From the definition of $\rho(m)$
 (recalling \eqref{countrho}), the congruence in \eqref{modlcm}
 is satisfied by exactly $\rho(\operatorname{lcm}(h, j))$ residue classes modulo $\operatorname{lcm}(h, j)$. 
 We are counting the integers in these classes that lie in
\begin{equation}\label{miscinterval} 
 \left[ \frac{\gcd(A, B)}{B}\left( \frac{x}{2a} - u_{0} \right), \frac{\gcd(A, B)}{B}\left( \frac{x}{a} - u_{0} \right) \right], 
\end{equation}
 noting that the interval in 
 \eqref{miscinterval} is of length
 $\frac{x \gcd(A, B)}{2 a B}$. 
 From a previous argument, the number of integers in each fixed residue class modulo $\operatorname{lcm}(h, j)$
 lying in the interval in \eqref{miscinterval} differs from 
 $\frac{x \gcd(A, B)}{2a B \operatorname{lcm}(h, j)}$ by at most $1$. 
 Since the $\rho(\operatorname{lcm}(h, j))$ residue classes are disjoint,
 by adding the number of integers belonging to each residue class,
 with only integers in the interval \eqref{miscinterval} being counted,
 we obtain exactly the inner sum over $k$ in \eqref{expandsquare}.

 Fix $h$ and $j$. Let $N_i$ denote the number of integers in \eqref{miscinterval} belonging to the $i^{\text{th}}$ of the 
 $\rho(\operatorname{lcm}(h, j))$ residue classes satisfying the congruence relation in \eqref{modlcm}. As above, we have established 
 that $$ \left| N_{i} - \frac{x \gcd(A, B)}{2 a B \, \operatorname{lcm}(h, j)} \right| \leq 1$$
 for each index $i$. 
 By construction, the inner sum in \eqref{expandsquare} over $k$ is precisely
 $ \sum_{i=1}^{\rho(\operatorname{lcm}(h, j))} N_{i}$. 
 An application of the triangle inequality then gives us that 
\begin{align*}
 & \left| \sum_{i=1}^{\rho(\operatorname{lcm}(h, j))} N_{i} 
 - \frac{x \gcd(A, B) \rho(\operatorname{lcm}(h, j))}{2 a B \, \operatorname{lcm}(h, j)} \right| \\ 
 & = \left| \sum_{i=1}^{\rho(\operatorname{lcm}(h, j))} \left( N_{i} 
 - \frac{x \gcd(A, B)}{2 a B \, \operatorname{lcm}(h, j)} \right) \right| \\ 
 & \leq \sum_{i=1}^{\rho(\operatorname{lcm}(h, j))}
 \left| N_{i} - \frac{x \gcd(A, B)}{2 a B \, \operatorname{lcm}(h, j)} \right| \\ 
 & \leq \sum_{i=1}^{\rho(\operatorname{lcm}(h, j))} 1 \\ 
 & = \rho(\operatorname{lcm}(h, j)).
\end{align*}
 So, we have showed that 
\begin{multline*}
 \frac{x\gcd(A,B)\,\rho(\operatorname{lcm}(h,j))}
 {2aB\,\operatorname{lcm}(h,j)}
-\rho(\operatorname{lcm}(h,j))
 \leq \sum_{i=1}^{\rho(\operatorname{lcm}(h,j))} N_i
 \leq \\ 
\frac{x\gcd(A,B)\,\rho(\operatorname{lcm}(h,j))}
 {2aB\,\operatorname{lcm}(h,j)}
+\rho(\operatorname{lcm}(h,j)). 
\end{multline*}
Multiplying the upper bound by $\lambda_h\lambda_j$ when
$\lambda_h\lambda_j \geq 0$, and the lower bound by
$\lambda_h\lambda_j$ when $\lambda_h\lambda_j < 0$, with the
inequality reversed in the latter case, gives
 $$ 
\lambda_h\lambda_j
\sum_{i=1}^{\rho(\operatorname{lcm}(h,j))} N_i
\leq {} 
\frac{x\gcd(A,B)}{2aB}\,
\lambda_h\lambda_j\,
\frac{\rho(\operatorname{lcm}(h,j))}
 {\operatorname{lcm}(h,j)} 
 + 
|\lambda_h\lambda_j|\,
\rho(\operatorname{lcm}(h,j)).
$$
Summing the last inequality over all pairs of positive integers
$h,j \leq z$ and using \eqref{expandsquare}, we obtain
\begin{align}
\begin{split}
 S \leq 
\frac{x\gcd(A,B)}{2aB} 
 & \sum_{1 \leq h,j \leq z}
\lambda_h\lambda_j\,
\frac{\rho(\operatorname{lcm}(h,j))}
 {\operatorname{lcm}(h,j)} + \\ 
 & \hspace{1.40in} \sum_{1 \leq h,j \leq z}
|\lambda_h\lambda_j|\,
\rho(\operatorname{lcm}(h,j)). 
\end{split}\label{Sdoublesum}
\end{align}

 Since $\frac{A}{\gcd(A, B)}$ and $\frac{B}{\gcd(A, B)}$ are relatively prime, for each prime $p$, at least one of the congruences among 
 \eqref{orcongruence1} and \eqref{orcongruence2}
 has a coefficient of $k$ that is invertible modulo $p$ 
 and therefore has a solution, so $\rho(p) \geq 1$.
 Now, for every positive squarefree integer $m$, we have
$\rho(m)=\prod_{p\mid m}\rho(p)>0$ by \eqref{furtherCRT}, and expanding the finite
product gives
\begin{align}
\begin{split}
\frac{m}{\rho(m)}
 &=
\prod_{p\mid m}\frac{p}{\rho(p)}
=
\prod_{p\mid m}
\left(1+\frac{p-\rho(p)}{\rho(p)}\right) \\ 
&=
\sum_{d\mid m}
\prod_{p\mid d}\frac{p-\rho(p)}{\rho(p)}. 
\end{split}\label{quotientCRT}
\end{align}
 For positive squarefree integers $h$ and $j$, the prime factorization formula for $\rho$, as given in \eqref{furtherCRT}, together with 
 the identity in \eqref{quotientCRT} applied to $m=\gcd(h,j)$, gives
\begin{align}
\begin{split}
\frac{\rho(\operatorname{lcm}(h,j))}
 {\operatorname{lcm}(h,j)}
&=
\frac{\rho(h)\rho(j)}{hj}\,
\frac{\gcd(h,j)}{\rho(\gcd(h,j))}
\\ 
&=
\frac{\rho(h)\rho(j)}{hj}
\sum_{\substack{d\mid h\\ d\mid j}}
\prod_{p\mid d}\frac{p-\rho(p)}{\rho(p)}.
\end{split}\label{lcmCRT}
\end{align}
Since $\lambda_h=0$ whenever $h$ is not squarefree,
only squarefree $h$ and $j$ contribute to the first double
sum in \eqref{Sdoublesum}. Applying \eqref{lcmCRT} to
these terms, we obtain
\begin{align*}
&\sum_{1\leq h,j\leq z}
\lambda_h\lambda_j
\frac{\rho(\operatorname{lcm}(h,j))}
{\operatorname{lcm}(h,j)}
\\
&\quad=
\sum_{\substack{1\leq h,j\leq z\\
h,j\ \mathrm{squarefree}}}
\frac{\lambda_h\rho(h)}{h}
\frac{\lambda_j\rho(j)}{j}
\sum_{\substack{d\mid h\\d\mid j}}
\prod_{p\mid d}\frac{p-\rho(p)}{\rho(p)}
\\
&\quad=
\sum_{\substack{1\leq d\leq z\\
d\ \mathrm{squarefree}}}
\left(\prod_{p\mid d}\frac{p-\rho(p)}{\rho(p)}\right)
\sum_{\substack{1\leq h,j\leq z\\
d\mid h,\ d\mid j}}
\frac{\lambda_h\rho(h)}{h}
\frac{\lambda_j\rho(j)}{j}.
\end{align*}
In the last equality, we have interchanged the finite sums.
Every common divisor $d$ of squarefree integers $h,j\leq z$
is itself squarefree and satisfies $d\leq z$.
We have also omitted the squarefree restrictions on the
inner indices $h$ and $j$, since any additional terms
vanish by the definition of the weights.

For each fixed $d$, the indices $h$ and $j$ range
independently over the same set of positive multiples of
$d$ not exceeding $z$. Thus the inner double sum factors
as
\begin{align*}
&\sum_{\substack{1\leq h,j\leq z\\
d\mid h,\ d\mid j}}
\frac{\lambda_h\rho(h)}{h}
\frac{\lambda_j\rho(j)}{j}
\\
&\quad=
\left(
\sum_{\substack{1\leq h\leq z\\d\mid h}}
\frac{\lambda_h\rho(h)}{h}
\right)
\left(
\sum_{\substack{1\leq j\leq z\\d\mid j}}
\frac{\lambda_j\rho(j)}{j}
\right)
\\
&\quad=
\left(
\sum_{\substack{1\leq h\leq z\\d\mid h}}
\frac{\lambda_h\rho(h)}{h}
\right)^2,
\end{align*}
where the last equality follows by renaming the summation
variable in the second factor. Consequently, we have that 
\begin{multline*} 
 \sum_{1\leq h,j\leq z}
\lambda_h\lambda_j\,
\frac{\rho(\operatorname{lcm}(h,j))}
 {\operatorname{lcm}(h,j)}
 = \\ 
\sum_{\substack{1\leq d\leq z\\ d\ \mathrm{squarefree}}}
\left(\prod_{p\mid d}\frac{p-\rho(p)}{\rho(p)}\right)
\left(
\sum_{\substack{1\leq h\leq z\\ d\mid h}}
\frac{\lambda_h\rho(h)}{h}
\right)^2. 
\end{multline*}

 For each positive squarefree integer $h\leq z$, we now specify 
\begin{equation}\label{fulllambda}
\lambda_h =
\frac{(-1)^{\omega(h)}}{\mathcal{G}(z)}
\left(\prod_{p\mid h}\frac{p}{p-\rho(p)}\right)
\sum_{\substack{
1\leq m\leq \frac{z}{h} \\
m\ \mathrm{squarefree}\\
\gcd(m,h)=1
}}
\prod_{p\mid m}\frac{\rho(p)}{p-\rho(p)},
\end{equation}
 while retaining $\lambda_h=0$ for all other positive integers $h$.
 Taking $h=1$ in \eqref{fulllambda}, 
 and using the relations $\omega(1)=0$ and $\gcd(m,1)=1$ together with 
 the empty-product convention,
we obtain
\[
\lambda_1
=
\frac{1}{\mathcal{G}(z)}
\sum_{\substack{1\leq m\leq z\\m\ \mathrm{squarefree}}}
\prod_{p\mid m}\frac{\rho(p)}{p-\rho(p)}
=
\frac{\mathcal{G}(z)}{\mathcal{G}(z)}
=
1.
\]
 For every positive squarefree integer $h\leq z$, the sum over $m$
 in the formula in \eqref{fulllambda} is at most $\mathcal{G}(z)$, since
$\frac{z}{h} \leq z$ and all the summands are nonnegative, so
\[
|\lambda_h|
\leq
\prod_{p\mid h}\frac{p}{p-\rho(p)}
\leq
\prod_{p\mid h}p
=
h,
\]
where the second inequality follows from $p-\rho(p)\geq 1$.

 Because $\rho(m)$ counts residue classes modulo $m$, we have $\rho(m)\leq m$ for every positive integer $m$, and hence, using
$\operatorname{lcm}(h,j)\leq hj$ and $|\lambda_h|\leq h$ also
when $\lambda_h=0$, we obtain
\[
\begin{aligned}
 \sum_{1\leq h,j\leq z}
|\lambda_h\lambda_j|\,
\rho(\operatorname{lcm}(h,j))
&\leq
\sum_{1\leq h,j\leq z} h^2j^2
\\[4pt]
&=
\left(\sum_{1\leq h\leq z} h^2\right)^2
\\[4pt]
&\leq
\left(\lfloor z\rfloor z^2\right)^2 \\ 
 & \leq z^6, 
\end{aligned}
\]
 so that \eqref{zx0005} then gives us the upper bound 
\begin{equation}\label{bound003}
 \sum_{1\leq h,j\leq z}
|\lambda_h\lambda_j|\,
\rho(\operatorname{lcm}(h,j))
 \leq x^{0.03} 
\end{equation}
 for the latter double sum in \eqref{Sdoublesum}. 

 For each positive squarefree integer $h$ not exceeding $ z$, multiplying \eqref{fulllambda} by $\frac{\rho(h)}{h}$, using the relation in 
 \eqref{furtherCRT} and then changing variables
from $m$ to $n=hm$, this gives us that 
\begin{align}
\begin{split}
\frac{\lambda_h\rho(h)}{h}
&= \frac{(-1)^{\omega(h)}}{\mathcal{G}(z)}
\left(\prod_{p\mid h}\frac{\rho(p)}{p-\rho(p)}\right)
\sum_{\substack{
1\leq m\leq \frac{z}{h} \\
m\ \mathrm{squarefree}\\
\gcd(m,h)=1
}}
\prod_{p\mid m}\frac{\rho(p)}{p-\rho(p)}
\\ 
&=
\frac{(-1)^{\omega(h)}}{\mathcal{G}(z)}
\sum_{\substack{
1\leq n\leq z\\
n\ \mathrm{squarefree}\\
h\mid n
}}
\prod_{p\mid n}\frac{\rho(p)}{p-\rho(p)}.
\end{split}\label{lambdarhoCRT}
\end{align}

 For each positive squarefree integer $d\leq z$, summing the preceding identity over squarefree multiples $h\leq z$ of $d$,
interchanging the finite sums, and recalling that $\lambda_h=0$
for nonsquarefree $h$, we obtain
\begin{equation}\label{sumsumprodsum}
 \sum_{\substack{1\leq h\leq z\\d\mid h}}
\frac{\lambda_h\rho(h)}{h}
 = 
\frac{1}{\mathcal{G}(z)}
\sum_{\substack{
1\leq n\leq z\\
n\ \mathrm{squarefree}\\
d\mid n
}}
\left(\prod_{p\mid n}\frac{\rho(p)}{p-\rho(p)}\right)
\sum_{\substack{h\mid n\\d\mid h}}
(-1)^{\omega(h)}.
\end{equation}

 For each positive squarefree integer $n$ divisible by $d$, the substitution $h=de$, with $e$ ranging over the positive divisors of 
 $\frac{n}{d}$, gives us that $\omega(h)=\omega(d)+\omega(e)$, since $\gcd(d,e)=1$, and hence
$$
\begin{aligned}
\sum_{\substack{h\mid n\\d\mid h}}
(-1)^{\omega(h)}
&=
(-1)^{\omega(d)}
\sum_{e\mid \frac{n}{d}}(-1)^{\omega(e)}
\\[4pt]
&=
(-1)^{\omega(d)}
\prod_{p\mid \frac{n}{d}}(1-1)
\\[4pt]
&=
\begin{cases}
(-1)^{\omega(d)},& n=d,\\
0,& n>d.
\end{cases}
\end{aligned}
$$
 For each positive squarefree integer $d\leq z$, substituting the above
 evaluation for $ \sum_{\substack{h\mid n\\d\mid h}}
(-1)^{\omega(h)} $ into the right-hand side of 
 \eqref{sumsumprodsum} leaves only the term corresponding to $n=d$, yielding 
$$
\sum_{\substack{1\leq h\leq z\\d\mid h}}
\frac{\lambda_h\rho(h)}{h}
=
\frac{(-1)^{\omega(d)}}{\mathcal{G}(z)}
\prod_{p\mid d}\frac{\rho(p)}{p-\rho(p)}.
$$
 Substituting this identity into the sum of squares for the first double sum in the bound for $S$
 displayed in \eqref{Sdoublesum}, and using
$(-1)^{2\omega(d)}=1$, we obtain
\[
\begin{aligned}
&\sum_{1\leq h,j\leq z}
\lambda_h\lambda_j\,
\frac{\rho(\operatorname{lcm}(h,j))}
 {\operatorname{lcm}(h,j)}
\\[4pt]
&\qquad=
\sum_{\substack{1\leq d\leq z\\d\ \mathrm{squarefree}}}
\left(\prod_{p\mid d}\frac{p-\rho(p)}{\rho(p)}\right)
\left(
\frac{(-1)^{\omega(d)}}{\mathcal{G}(z)}
\prod_{p\mid d}\frac{\rho(p)}{p-\rho(p)}
\right)^2
\\[4pt]
&\qquad=
\frac{1}{\mathcal{G}(z)^2}
\sum_{\substack{1\leq d\leq z\\d\ \mathrm{squarefree}}}
\prod_{p\mid d}\frac{\rho(p)}{p-\rho(p)}
\\[4pt]
&\qquad 
 =
\frac{1}{\mathcal{G}(z)}.
\end{aligned}
\]
 Substituting $\frac{1}{\mathcal{G}(z)}$ for the first double sum
 in \eqref{Sdoublesum}, and using the upper bound of $x^{0.03}$ in \eqref{bound003} 
for the latter double sum, we obtain
\begin{equation}\label{Sboundx003}
 S \leq \frac{x\gcd(A,B)}{2aB\,\mathcal{G}(z)} +x^{0.03}. 
\end{equation}
 Since $\rho(p)=2$ whenever $p\nmid 2AB|a-b|$, omitting the nonnegative terms with $\gcd(h,2AB|a-b|)>1$ from
$\mathcal{G}(z)$ and using the relation $\frac{2}{p-2} \geq \frac{2}{p}$ for $p>2$ gives
$$
\begin{aligned}
\mathcal{G}(z)
&\geq
\sum_{\substack{
1\leq h\leq z\\
h\ \mathrm{squarefree}\\
\gcd(h,2AB|a-b|)=1
}}
\prod_{p\mid h}\frac{2}{p-2}
\\[4pt]
&\geq
\sum_{\substack{
1\leq h\leq z\\
h\ \mathrm{squarefree}\\
\gcd(h,2AB|a-b|)=1
}}
\frac{2^{\omega(h)}}{h}.
\end{aligned}
$$

 Writing each positive squarefree integer $h$ uniquely as $h=dm$, where $d$ is the product of the prime divisors of $h$ that divide
 $2AB|a-b|$ and $m = \frac{h}{d}$, and using
 $\omega(h)=\omega(d)+\omega(m)$, we obtain
$$
 \sum_{\substack{1\leq h\leq z\\h\ \mathrm{squarefree}}}
\frac{2^{\omega(h)}}{h}
 =
\sum_{\substack{
d\mid 2AB|a-b|\\
d\ \mathrm{squarefree}
}}
\frac{2^{\omega(d)}}{d}
\sum_{\substack{
1\leq m\leq \frac{z}{d} \\
m\ \mathrm{squarefree}\\
\gcd(m,2AB|a-b|)=1
}}
\frac{2^{\omega(m)}}{m}.
$$

 Since $\frac{z}{d} \leq z$ and all summands are nonnegative, replacing $m\leq \frac{z}{d}$ by $m\leq z$ in each inner sum and using 
 the expansion of $\mathcal{B}_{1,2}$ over squarefree divisors gives
\begin{align*}
 \sum_{\substack{1\leq h\leq z\\h\ \mathrm{squarefree}}}
\frac{2^{\omega(h)}}{h} 
 & \leq
\left(
\sum_{\substack{d\mid 2AB|a-b|\\d\ \mathrm{squarefree}}}
\frac{2^{\omega(d)}}{d}
\right)
\sum_{\substack{
1\leq m\leq z\\
m\ \mathrm{squarefree}\\
\gcd(m,2AB|a-b|)=1
}}
\frac{2^{\omega(m)}}{m} \\ 
& =
\mathcal{B}_{1,2}(2AB|a-b|)
\sum_{\substack{
1\leq m\leq z\\
m\ \mathrm{squarefree}\\
\gcd(m,2AB|a-b|)=1
}}
\frac{2^{\omega(m)}}{m}.
\end{align*}
 Since the sum over $m$ on the right is at most $\mathcal{G}(z)$, dividing by the positive quantity
 $\mathcal{B}_{1,2}(2AB|a-b|)$ gives
\begin{equation}\label{compareGandB}
\mathcal{G}(z)
\geq
\frac{1}{\mathcal{B}_{1,2}(2AB|a-b|)}
\sum_{\substack{1\leq h\leq z\\h\ \mathrm{squarefree}}}
\frac{2^{\omega(h)}}{h}.
\end{equation}
 Since every positive squarefree integer $h$ has exactly
$2^{\omega(h)}$ ordered factorizations $h=mn$ into positive
integers, and restricting to odd $m,n\leq\sqrt z$ only removes
nonnegative terms, we have
$$
\begin{aligned}
\sum_{\substack{1\leq h\leq z\\h\ \mathrm{squarefree}}}
\frac{2^{\omega(h)}}{h}
&=
\sum_{\substack{
m,n\geq 1\\
mn\leq z\\
mn\ \mathrm{squarefree}
}}
\frac{1}{mn}
\\[4pt]
&\geq
\sum_{\substack{
1\leq m,n\leq\sqrt z\\
m,n\ \mathrm{odd}\\
mn\ \mathrm{squarefree}
}}
\frac{1}{mn}.
\end{aligned}
$$ Separating the ordered pairs of odd positive integers $m,n\leq\sqrt z$ according to whether $mn$ is squarefree,
and factoring the sum over all such pairs, gives
$$
 \sum_{\substack{
1\leq m,n\leq\sqrt z\\
m,n\ \mathrm{odd}\\
mn\ \mathrm{squarefree}
}}
 \frac{1}{mn}
 =
\left(
\sum_{\substack{1\leq m\leq\sqrt z\\m\ \mathrm{odd}}}
\frac{1}{m}
\right)^2
-
\sum_{\substack{
1\leq m,n\leq\sqrt z\\
m,n\ \mathrm{odd}\\
mn\ \text{not squarefree}
}}
\frac{1}{mn}.
$$ Now, if $m$ and $n$ are odd positive integers and $mn$ is not squarefree, then some odd prime $p$ satisfies $p^2\mid mn$, 
 and considering the exponents of $p$ in $m$ and $n$ shows
that at least one proposition among 
$ p^2\mid m$, $p^2\mid n$, and 
 $ \bigl(p\mid m\ \text{and}\ p\mid n\bigr)$ holds.

 For each odd prime $p$, substituting $m=p^2k$ and using the relation $\frac{\sqrt z}{p^2} \leq\sqrt z$ gives $$
\begin{aligned}
\sum_{\substack{
1\leq m,n\leq\sqrt z\\
m,n\ \mathrm{odd}\\
p^2\mid m
}}
\frac{1}{mn}
&=
\frac{1}{p^2}
\left(
\sum_{\substack{1\leq k\leq \frac{\sqrt{z}}{p^2} \\k\ \mathrm{odd}}}
\frac{1}{k}
\right)
\left(
\sum_{\substack{1\leq n\leq\sqrt z\\n\ \mathrm{odd}}}
\frac{1}{n}
\right)
\\[4pt]
&\leq
\frac{1}{p^2}
\left(
\sum_{\substack{1\leq m\leq\sqrt z\\m\ \mathrm{odd}}}
\frac{1}{m}
\right)^2,
\end{aligned}
$$
 and interchanging $m$ and $n$ gives the same upper bound when
$p^2\mid m$ is replaced by $p^2\mid n$.

For each odd prime $p$, separating the sums over $m$ and $n$,
substituting $m=pk$, and using $\frac{\sqrt z}{p} \leq\sqrt z$, we obtain
\[
\begin{aligned}
\sum_{\substack{
1\leq m,n\leq\sqrt z\\
m,n\ \mathrm{odd}\\
p\mid m,\ p\mid n
}}
\frac{1}{mn}
&=
\left(
\sum_{\substack{
1\leq m\leq\sqrt z\\
m\ \mathrm{odd}\\
p\mid m
}}
\frac{1}{m}
\right)^2
\\[4pt]
&=
\frac{1}{p^2}
\left(
\sum_{\substack{1\leq k\leq \frac{\sqrt z}{p} \\k\ \mathrm{odd}}}
\frac{1}{k}
\right)^2
\\[4pt]
&\leq
\frac{1}{p^2}
\left(
\sum_{\substack{1\leq m\leq\sqrt z\\m\ \mathrm{odd}}}
\frac{1}{m}
\right)^2.
\end{aligned}
\]

 Since every ordered pair of odd positive integers $m,n\leq\sqrt z$ with $mn$ not squarefree satisfies $p^2\mid m$, $p^2\mid n$, or 
 both $p\mid m$ and $p\mid n$ for some prime $2<p\leq\sqrt z$, summing the three bounds over these primes, with repeated 
 counting contributing only additional nonnegative terms, gives
$$
\sum_{\substack{
1\leq m,n\leq\sqrt z\\
m,n\ \mathrm{odd}\\
mn\ \text{not squarefree}
}}
\frac{1}{mn}
\leq
3\left(\sum_{2<p\leq\sqrt z}\frac{1}{p^2}\right)
\left(
\sum_{\substack{1\leq m\leq\sqrt z\\m\ \mathrm{odd}}}
\frac{1}{m}
\right)^2.
$$
 Euler's evaluation for $1 + \frac{1}{2^2} + \frac{1}{3^2} + \cdots = \frac{\pi^2}{6}$ 
 allows us to obtain the bound 
\begin{equation}\label{boundBasel}
 \sum_{2<p\leq\sqrt z}\frac{1}{p^2}
 \leq 
\frac{1}{4}.
\end{equation}
 Using \eqref{boundBasel} and the
bound for the sum with $mn$ not squarefree gives
\[
\begin{aligned}
\sum_{\substack{
1\leq m,n\leq\sqrt z\\
m,n\ \mathrm{odd}\\
mn\ \mathrm{squarefree}
}}
\frac{1}{mn}
&\geq
\left(1-\frac{3}{4}\right)
\left(
\sum_{\substack{1\leq m\leq\sqrt z\\m\ \mathrm{odd}}}
\frac{1}{m}
\right)^2
\\[4pt]
&=
\frac{1}{4}
\left(
\sum_{\substack{1\leq m\leq\sqrt z\\m\ \mathrm{odd}}}
\frac{1}{m}
\right)^2.
\end{aligned}
\]
 A standard/classical estimate for skew/odd harmonic numbers
 (which are of the form $1 + \frac{1}{3} + \cdots + \frac{1}{2n+1}$ for nonnegative integers $n$) 
 gives us that 
\begin{equation}\label{classicalskew}
 \sum_{\substack{1\leq m\leq\sqrt z\\m\ \mathrm{odd}}}
\frac{1}{m}
 \geq \frac{1}{4} \log z. 
\end{equation} 
 Combining the lower bound obtained from odd factor pairs with
the estimate in \eqref{classicalskew}, and noting that $z\geq 1$ implies
$\log z\geq 0$, we obtain
\[
\begin{aligned}
\sum_{\substack{1\leq h\leq z\\h\ \mathrm{squarefree}}}
\frac{2^{\omega(h)}}{h}
&\geq
\frac{1}{4}
\left(
\sum_{\substack{1\leq m\leq\sqrt z\\m\ \mathrm{odd}}}
\frac{1}{m}
\right)^2
\\[4pt]
&\geq
\frac{1}{4}\left(\frac{1}{4}\log z\right)^2
=
\frac{(\log z)^2}{64}.
\end{aligned}
\]
 Combining this estimate with the lower bound for $\mathcal{G}(z)$ displayed in \eqref{compareGandB}, we obtain
 \[
\begin{aligned}
\mathcal{G}(z)
&\geq
\frac{1}{\mathcal{B}_{1,2}(2AB|a-b|)}
\sum_{\substack{1\leq h\leq z\\h\ \mathrm{squarefree}}}
\frac{2^{\omega(h)}}{h}
\\[4pt]
&\geq
\frac{(\log z)^2}
 {64\,\mathcal{B}_{1,2}(2AB|a-b|)}.
\end{aligned}
\]
 Since $z=x^{0.005}>1$, substituting this lower bound for
$\mathcal{G}(z)$ into the upper bound for $S$ shown in \eqref{Sboundx003} gives
\[
\begin{aligned}
S
&\leq
\frac{x\gcd(A,B)}{2aB\,\mathcal{G}(z)}
+x^{0.03}
\\[4pt]
&\leq
\frac{x\gcd(A,B)}{2aB}\,
\frac{64\,\mathcal{B}_{1,2}(2AB|a-b|)}
 {(\log z)^2}
+x^{0.03}
\\[4pt]
&=
\frac{32x\gcd(A,B)}
 {aB(\log z)^2}\,
\mathcal{B}_{1,2}(2AB|a-b|)
+x^{0.03}.
\end{aligned}
\]
 Applying Lemma~\ref{lemcalBineq} with $\beta=1$, $K=2$,
 $M=2$, and $N=AB|a-b|$, and using 
 $\mathcal{B}_{1,2}(2)=2$, we obtain 
\begin{equation}\label{applysubmult} 
 \mathcal{B}_{1,2}(2AB|a-b|)
\leq
\mathcal{B}_{1,2}(2)\,
\mathcal{B}_{1,2}(AB|a-b|)
=
2\,\mathcal{B}_{1,2}(AB|a-b|).
\end{equation} 
Using \eqref{applysubmult}, together with
$B\geq \frac{b}{2}$ and $\log z=\frac{\log x}{200}$,
we obtain
 \[
\begin{aligned}
S
&\leq
\frac{64x\gcd(A,B)}
 {aB(\log z)^2}\,
\mathcal{B}_{1,2}(AB|a-b|)
+x^{0.03}
\\[4pt]
&\leq
\frac{128x\gcd(A,B)}
 {ab(\log z)^2}\,
\mathcal{B}_{1,2}(AB|a-b|)
+x^{0.03}
\\[4pt]
&=
\frac{5\,120\,000\,x\gcd(A,B)}
 {ab(\log x)^2}\,
\mathcal{B}_{1,2}(AB|a-b|)
+x^{0.03}.
\end{aligned}
\]
 
 Since the number of ordered prime pairs satisfying the lemma's conditions is at most $S$, the required estimate follows with 
 $C=5\,120\,000$ and $x_0=8$, completing the proof.
\end{proof}

\begin{lemma}\label{lemlineargcd}
 Let $U$ and $c$ be positive integers satisfying $\gcd(U,c)=1$, and let $R \geq 2$ be real. Then
$$
 \sum_{\substack{1 \leq m,n \leq R\\m \neq n}}
 \gcd(Um+c,Un+c)
 \ll R^2 \log R,
$$ where the sum is over ordered pairs of positive integers,
and the implied constant is absolute. In particular, the
bound is uniform in $U$, $c$, and $R$.
\end{lemma}

\begin{proof}
For every positive integer $m$, we have
\begin{equation}\label{displaycoprime}
 \gcd(U,Um+c)=\gcd(U,c)=1. 
\end{equation} 
 Suppose that $m$ and $n$ are distinct positive integers
not exceeding $R$, and that a positive integer $h$ divides
both $Um+c$ and $Un+c$. Since $h \mid Um+c$, the 
 coprimality relations shown in \eqref{displaycoprime} 
 imply that $\gcd(h,U)=1$. Also,
subtracting the two integers divisible by $h$ gives
$
 h \mid \bigl((Um+c)-(Un+c)\bigr)=U(m-n)$. 
Because $\gcd(h,U)=1$, it follows that $h \mid m-n$.
The assumption $m \neq n$ therefore gives
$
 1 \leq h \leq |m-n| \leq R-1 < R$. 
Thus every positive common divisor of $Um+c$ and $Un+c$
is less than $R$ and is relatively prime to $U$.

For each positive integer $h<R$ satisfying $\gcd(h,U)=1$,
define
\[
 N_h
 =
 \left|
 \left\{
 m \in \mathbb Z:
 1 \leq m \leq R,\quad h \mid Um+c
 \right\}
 \right|.
\]
We first bound $N_h(N_h-1)$.

If $N_h \geq 1$, write the integers counted by $N_h$ in
increasing order as
$
 m_1<m_2<\cdots<m_{N_h}$. 
For $1 \leq i<N_h$, we have
$
 h \mid
 \bigl((Um_{i+1}+c)-(Um_i+c)\bigr)$. 
Since $\gcd(h,U)=1$, this implies that
$h \mid (m_{i+1}-m_i)$. Each difference $m_{i+1}-m_i$ is positive, so $ m_{i+1}-m_i \geq h$. Adding these inequalities yields
$
 (N_h-1)h
 \leq m_{N_h}-m_1
 \leq R-1
 < R$. 
This inequality also holds when $N_h=1$, since its
left-hand side is then zero. Consequently, each of the relations among 
$
 N_h-1 \leq \frac{R}{h}$ and 
$ N_h \leq \frac{R}{h}+1$ holds. Both $N_h$ and $N_h-1$ are nonnegative when $N_h \geq 1$,
so we may multiply these bounds to obtain
\begin{equation}\label{extendNh0}
 N_h(N_h-1)
 \leq
 \left(\frac{R}{h}+1\right)\frac{R}{h}
 =
 \frac{R^2}{h^2}+\frac{R}{h}.
\end{equation}
 If $N_h=0$, the same bound in \eqref{extendNh0} holds, because its left-hand
 side vanishes. We have therefore proved that
\begin{equation}\label{eqlineargcdcount}
 N_h(N_h-1)
 \leq
 \frac{R^2}{h^2}+\frac{R}{h}
\end{equation}
for every positive integer $h<R$ with $\gcd(h,U)=1$.

For distinct positive integers $m,n \leq R$, the greatest
common divisor of $Um+c$ and $Un+c$ is itself one of
their positive common divisors. Since all positive common
divisors are less than $R$ and relatively prime to $U$,
we have
\[
 \gcd(Um+c,Un+c)
 \leq
 \sum_{\substack{1 \leq h<R\\
 \gcd(h,U)=1\\
 h \mid Um+c\\
 h \mid Un+c}}
 h.
\]
Indeed, the sum on the right includes the term
$h=\gcd(Um+c,Un+c)$, and every other term is positive.

For each fixed $h$, the number of ordered pairs of
distinct integers counted by $N_h$ is $N_h(N_h-1)$.
Therefore, summing the preceding inequality over $m,n$
and interchanging the finite sums gives
\begin{align*}
 &\sum_{\substack{1 \leq m,n \leq R\\m \neq n}}
 \gcd(Um+c,Un+c)
 \\
 &\qquad\leq
 \sum_{\substack{1 \leq m,n \leq R\\m \neq n}}
 \sum_{\substack{1 \leq h<R\\
 \gcd(h,U)=1\\
 h \mid Um+c\\
 h \mid Un+c}}
 h
 \\
 &\qquad=
 \sum_{\substack{1 \leq h<R\\\gcd(h,U)=1}}
 h
 \sum_{\substack{1 \leq m,n \leq R\\
 m \neq n\\
 h \mid Um+c\\
 h \mid Un+c}}
 1
 \\
 &\qquad=
 \sum_{\substack{1 \leq h<R\\\gcd(h,U)=1}}
 h\,N_h(N_h-1).
\end{align*}
Applying \eqref{eqlineargcdcount}, and then including the
nonnegative terms corresponding to integers $h$ with
$\gcd(h,U)>1$, we obtain
\begin{align*}
 \sum_{\substack{1 \leq m,n \leq R\\m \neq n}}
 \gcd(Um+c,Un+c)
 &\leq
 \sum_{\substack{1 \leq h<R\\\gcd(h,U)=1}}
 h\left(\frac{R^2}{h^2}+\frac{R}{h}\right)
 \\
 &\leq
 \sum_{1 \leq h<R}
 h\left(\frac{R^2}{h^2}+\frac{R}{h}\right)
 \\
 &=
 R^2 \sum_{1 \leq h<R}\frac{1}{h}
 +
 R \sum_{1 \leq h<R}1.
\end{align*}
 Classical bounds for harmonic numbers 
 give us that 
 $ \sum_{1 \leq h<R}\frac{1}{h}
 \leq
 1+\log R$. 
Also, we see that 
$
 \sum_{1 \leq h<R}1
 \leq \lfloor R\rfloor
 \leq R$. 
 Combining these estimates gives 
\begin{equation}\label{floorharmonic}
 \sum_{\substack{1 \leq m,n \leq R\\m \neq n}}
 \gcd(Um+c,Un+c)
 \leq
 R^2(1+\log R)+R^2
 =
 R^2(2+\log R).
\end{equation}
Finally, since $R \geq 2$, we find that 
$
 2+\log R
 \leq
 \left(1+\frac{2}{\log 2}\right)\log R$. 
 Thus, from \eqref{floorharmonic}, we find that 
\[
 \sum_{\substack{1 \leq m,n \leq R\\m \neq n}}
 \gcd(Um+c,Un+c)
 \leq
 \left(1+\frac{2}{\log 2}\right)R^2\log R.
\]
 The constant on the right is independent of $U$, $c$, and $R$, as required.
\end{proof}

\begin{theorem}\label{maintheorem}
There exist absolute constants $c_0>0$ and $x_0\geq1$ such that $ V(x)\geq c_0x$
 for every real $x\geq x_0$. 
 Consequently, the relation 
\begin{equation}\label{strongest}
 \liminf_{x\to\infty}\frac{V(x)}x>0 
\end{equation}
 holds, and, in particular, the relation $V(x)\neq o(x)$ holds.
\end{theorem}

\begin{proof}
 Let $x$ be a sufficiently large real number, and put
\begin{equation}\label{x1000th}
 Y=x^{\frac{1}{1000}}.
\end{equation}
We take $x$ large enough that $Y\geq2$.
Every further lower bound on $x$ imposed below will depend only
on fixed numerical constants and on the uniform constants in
the preceding lemmas.

 Define $$ \mathcal E =\{e\in\mathbb Z:Y<e\leq Y^2,\ e\equiv2\pmod4\}, $$ and define the two sets of primes
\[
 \mathcal T=\{p:Y^{10}<p\leq Y^{11}\},
 \qquad
 \mathcal R=\{p:Y^{40}<p\leq Y^{41}\}.
\]
Set
\[
 \mathcal C=\{et:e\in\mathcal E,\ t\in\mathcal T\},
 \qquad
 \mathcal A=\{cr:c\in\mathcal C,\ r\in\mathcal R\},
\]
and
\begin{equation}\label{definecalNx}
 \mathcal N_x
 =\left\{ap:a\in\mathcal A,\quad p\ \text{prime},\quad
 \frac{x}{2a}<p\leq\frac xa\right\}.
\end{equation}
 These definitions imply
\begin{equation}\label{mainfactorranges}
 Y^{11}<c\leq Y^{13}\quad(c\in\mathcal C),
 \qquad
 Y^{51}<a\leq Y^{54}\quad(a\in\mathcal A).
\end{equation}
In particular, every prime divisor of $e\in\mathcal E$ is
smaller than every $t\in\mathcal T$, and every such $t$ is
smaller than every $r\in\mathcal R$.
 Moreover, for $ap\in\mathcal N_x$, we have that $$
 p>\frac{x}{2a}\geq\frac12Y^{946}>Y^{41}.
$$
Thus, in the factorization $ap=etrp$, the primes $t,r,p$
are, in increasing order, the three largest prime divisors,
and each occurs with exponent one.
 The fundamental theorem of arithmetic therefore implies that every
 factorization among 
\begin{align} 
 c & = e t, \label{cequalset} \\ 
 a & = c r = e t r, \quad \text{and} \label{aequalscr} \\ 
 n & = a p = e t r p \label{nequalsap}
\end{align}
 in the
 sets above is unique.
By the uniqueness asserted for the factorization $a=cr$
in \eqref{aequalscr}, we have
\begin{equation}\label{mainuniqueab}
 cr=ds
 \quad\Longleftrightarrow\quad
 c=d\ \text{and}\ r=s
 \qquad(c,d\in\mathcal C,\ r,s\in\mathcal R).
\end{equation}
 Also,
\begin{equation}\label{maininterval}
 \mathcal N_x\subseteq\Big(\frac{x}{2},x\Big]\cap\mathbb Z.
\end{equation}
Since $e\equiv2\pmod4$ and $t,r$ are odd, we have
\begin{equation}\label{maintwodivides}
 2\parallel c\quad(c\in\mathcal C),
 \qquad
 2\parallel a\quad(a\in\mathcal A).
\end{equation}

 By the prime number theorem, the asymptotic equivalence $
 \pi(v) - \pi\left(\frac{v}{2}\right)\sim\frac{v}{2\log v}
$
 holds. So, there exists an absolute $v_0$ such that
\begin{equation}\label{explicitPNT}
 \pi(v) - \pi\left(\frac{v}{2}\right)\geq\frac{v}{4\log v}
\end{equation}
 for all $v\geq v_0$. 
 Now, for every $a\in\mathcal A$, we have
$ x^{0.946}\leq\frac xa\leq x$. 
 Consequently, for all sufficiently large $x$, the 
 inequality in \eqref{explicitPNT} applies simultaneously to each $v = \frac{x}{a}$, giving
\begin{equation}\label{mainouterprimecount}
 \pi\left(\frac{x}{a}\right)-\pi\left( \frac{x}{2a} \right)
 \geq\frac{x}{4a\log\left(\frac{x}{a}\right)}
 \geq\frac{x}{4a\log x}.
\end{equation}

For $v\geq1$, the number of positive integers $e\leq v$
with $e\equiv2\pmod4$ is
$
 \left\lfloor\frac{v+2}{4}\right\rfloor
 =\frac v4+O(1)$. 
Partial summation therefore gives
\begin{equation}\label{harmonice}
 \sum_{e\in\mathcal E}\frac1e = 
 \frac{\left\lfloor \frac{Y^2+2}{4} \right\rfloor}{Y^2}
 -\frac{\left\lfloor \frac{Y+2}{4} \right\rfloor}{Y}
 +\int_Y^{Y^2}
 \frac{ \left\lfloor \frac{v+2}{4} \right\rfloor}{v^2}\,dv 
 =\frac14\log Y+O(Y^{-1}). 
\end{equation}
 For either $(\eta,\zeta)=(10,11)$ or $(\eta,\zeta)=(40,41)$,
partial summation also gives
\[
 \sum_{Y^\eta<p\leq Y^\zeta}\frac1p
 =
 \frac{\pi(Y^\zeta)}{Y^\zeta}
 -\frac{\pi(Y^\eta)}{Y^\eta}
 +\int_{Y^\eta}^{Y^\zeta}\frac{\pi(v)}{v^2}\,dv.
\]
The two endpoint terms tend to zero by the prime number
theorem. The error in replacing the integral by
$\int_{Y^\eta}^{Y^\zeta} \frac{dv}{v\log v}$ is at most
\[
 \left(
 \sup_{v\geq Y^\eta}
 \left|\frac{\pi(v)\log v}{v}-1\right|
 \right)
 \int_{Y^\eta}^{Y^\zeta}\frac{dv}{v\log v}
 =o(1),
\]
 since the supremum tends to zero and the integral equals $\log\left( \frac{\zeta}{\eta} \right)$. It follows that
\begin{equation}\label{mainprimemass}
 \sum_{t\in\mathcal T}\frac1t
 =\log\left( \frac{11}{10} \right)+o(1),
 \qquad
 \sum_{r\in\mathcal R}\frac1r
 =\log\left( \frac{41}{40} \right)+o(1).
\end{equation}

The uniqueness asserted for the factorization $n=ap$ in \eqref{nequalsap} gives us that $$ |\mathcal N_x|
 =
 \sum_{a\in\mathcal A}
 \ \left(\pi\left( \frac{x}{a} \right) - 
 \pi\left( \frac{x}{2a} \right) \right), 
$$ and the
 uniqueness asserted for the factorization $a=etr$ in
\eqref{aequalscr} gives us that 
 \[
 \sum_{a\in\mathcal A}\frac1a
 =
 \left(\sum_{e\in\mathcal E}\frac1e\right)
 \left(\sum_{t\in\mathcal T}\frac1t\right)
 \left(\sum_{r\in\mathcal R}\frac1r\right).
\]
 Combining \eqref{mainouterprimecount}, \eqref{harmonice}, and \eqref{mainprimemass},
and using $\log Y = \frac{\log x}{1000}$, gives
\begin{equation}\label{mainmass}
 |\mathcal N_x|
 \gg\frac{x\log Y}{\log x} \gg x.
\end{equation}

 For distinct $a,b\in\mathcal A$, define $ G(a,b)=\gcd(w(a),w(b),|a-b|)$ and $ D(a,b)=w(a)w(b)|a-b|$. 
 Observe that both $G(a,b)$ and $D(a,b)$ 
 are positive integers.
 Fix distinct $c,d\in\mathcal C$, and let 
 $ U=w(c)$ and $ W=w(d)$. 
 Since $r>c$ and $s>d$ whenever $r, s \in\mathcal R$,
 Lemma \ref{lemmawap} gives that $ w(cr)=Ur+c$ and that $ w(ds)=Ws+d$. By \eqref{mainuniqueab},  
  the relation   $cr\ne ds$   holds  for all such $r,s$.

 Suppose that $h\mid G(cr,ds)$. Applying Lemma~\ref{gcdmwm} with $n=cr$ and $n=ds$, we obtain
$\gcd(cr,w(cr))=\gcd(ds,w(ds))=1$.
Since $h\mid w(cr)$ and $h\mid w(ds)$, it follows that
$\gcd(h,cr)=\gcd(h,ds)=1$. 
 We thus find that $h$ and $cdrs$ are relatively prime. 
 If a prime were to divide both $h$ and $U$, then the relation
  $h  \mid (Ur+c)$     would force that same prime to divide $c$,
contrary to $\gcd(h,  c) = 1$. So, we can deduce that $h$ and $U$ are coprime. 
 The corresponding argument gives  us that  $\gcd(h,W)=1$, so
\begin{equation}\label{mainhcoprime}
 \gcd(h,cdUW)=1.
\end{equation}
The congruences
\[
 Ur\equiv-c\pmod h,
 \qquad
 Ws\equiv-d\pmod h
\]
therefore each specify one reduced residue class modulo $h$.
Here the residue classes are reduced because both $c$ and $d$
are relatively prime to $h$.
Moreover, $h\mid cr-ds$. Multiplying
$cr\equiv ds\pmod h$ by $UW$ and substituting the preceding
two congruences gives
\[
 -c^2W\equiv-d^2U\pmod h.
\]
 As a consequence, each of the relations among 
\begin{equation}\label{mainHdefinition}
 h\mid H(c,d) \quad \text{and} \quad 
 H(c,d)=c^2w(d)-d^2w(c) 
\end{equation}
 holds. 

 The integer $H(c,d)$ is nonzero. Indeed, if it were zero, then $c^2w(d)=d^2w(c)$. Applying Lemma~\ref{gcdmwm}
with $n=c$ and $n=d$, we obtain
$\gcd(c,w(c))=\gcd(d,w(d))=1$.
Consequently, the equalities among 
\begin{equation}\label{miscgcdeqs} 
 \gcd(c^2,w(c))=\gcd(d^2,w(d))=1
\end{equation}
 hold, and from \eqref{miscgcdeqs}, we find that 
 the equality $c^2w(d)=d^2w(c)$ would imply that 
 $c^2\mid d^2$ and that $d^2\mid c^2$.
Hence $c^2=d^2$, and positivity would imply $c=d$,
contrary to our choice of $c,d$.
 Also, applying the upper bound in Lemma \ref{lemboundw} with $n = c$ and $n = d$, 
 and using $c, d \leq Y^{13}$, we obtain 
\begin{align}
 0<|H(c,d)|
 &\leq c^2w(d)+d^2w(c)
 \notag\\
 &\leq\frac{26}{\log2}Y^{39}\log Y
 <Y^{39.5}
 \label{mainHsize}
\end{align}
for all sufficiently large $Y$, uniformly in $c,d$.

The sum over the positive divisors of $G(cr,ds)$ includes
the term corresponding to $h=G(cr,ds)$, so
\[
 G(cr,ds)^{\frac{3}{2}}
 \leq\sum_{h\mid G(cr,ds)}h^{\frac{3}{2}}.
\]
Using \eqref{mainhcoprime} and \eqref{mainHdefinition},
interchanging finite sums, and then dropping the condition
$h\mid cr-ds$, we obtain
\begin{align*}
 \sum_{r,s\in\mathcal R}\frac{G(cr,ds)^{\frac{3}{2}}}{rs}
 &\leq
 \sum_{\substack{h\mid |H(c,d)|\\\gcd(h,cdUW)=1}}
 h^{\frac{3}{2}}
 \sum_{\substack{r,s\in\mathcal R\\
 h\mid Ur+c\\h\mid Ws+d\\h\mid cr-ds}}
 \frac1{rs}
 \\
 &\leq
 \sum_{\substack{h\mid |H(c,d)|\\\gcd(h,cdUW)=1}}
 h^{\frac{3}{2}}
 \left(\sum_{\substack{r\in\mathcal R\\h\mid Ur+c}}\frac1r\right)
 \left(\sum_{\substack{s\in\mathcal R\\h\mid Ws+d}}\frac1s\right).
\end{align*}
Each inner sum is over one reduced residue class modulo $h$.
 We then choose its representative to be a positive integer, as
required in Lemma~\ref{lemreciprocal}.
By \eqref{mainHsize}, Lemma~\ref{lemreciprocal} applies with
$ \alpha=40$, $\delta=41$, and $\gamma=39.5$. 
 Its constant is independent of $c,d,U,W,h$.
 After applying the Lemma~\ref{lemreciprocal}, we may include the remaining
 divisors of $|H(c,d)|$, since all terms are nonnegative.
This gives
\begin{equation}\label{mainconditionaldivisors}
 \sum_{r,s\in\mathcal R}\frac{G(cr,ds)^{\frac{3}{2}}}{rs}
 \ll
 \sum_{h\mid |H(c,d)|}\frac{h^{\frac{3}{2}}}{\varphi(h)^2}.
\end{equation}

For every positive integer $N$, Euler's formula for
$\varphi$ and the fundamental theorem of arithmetic give
\begin{align}
\begin{split}
 \sum_{h\mid N}\frac{h^{\frac{3}{2}}}{\varphi(h)^2}
 &=
 \prod_{p^\nu\parallel N}
 \left(
 1+\frac1{(1-p^{-1})^2}
 \sum_{j=1}^{\nu}p^{-\frac{j}{2}}
 \right)
 \\
 &\leq
 \prod_{p\mid N}
 \left(
 1+\frac{p^{-\frac{1}{2}}}
 {(1-p^{-1})^2(1-p^{-\frac{1}{2}})}
 \right).
\end{split}\label{phisumtoprod}
\end{align}
For every $p\geq2$, the relations among 
$
 \frac1{(1-p^{-1})^2(1-p^{-\frac{1}{2}})}
 \leq\frac4{1-2^{-\frac{1}{2}}}
 =8+4\sqrt2<14$ hold. 
 This and \eqref{phisumtoprod} together give us that 
\begin{equation}\label{phisumtoB} 
 \sum_{h\mid N}\frac{h^{\frac{3}{2}}}{\varphi(h)^2}
 \leq\mathcal B_{\frac{1}{2},14}(N) 
\end{equation}
 for $N \geq 2$, and we may check that 
 \eqref{phisumtoB} also holds for 
 $N=1$ (in which case both sides of \eqref{phisumtoB} equal $1$). 
 We then see that \eqref{mainconditionaldivisors} and \eqref{phisumtoB} together imply 
\begin{equation}\label{mainconditionalG}
 \sum_{r,s\in\mathcal R}\frac{G(cr,ds)^{\frac{3}{2}}}{rs} \ll \mathcal B_{\frac{1}{2},14}(|H(c,d)|).
\end{equation}

 Our goal, at this point, is to estimate a sum involving the right-hand side of \eqref{mainconditionalG}
 for distinct $c, d \in \mathcal{C}$. As an intermediary step, 
 we begin by bounding an associated sum over elements $t \in \mathcal{T}$ 
 such that $et \neq fu$. Fix $e,f\in\mathcal E$ and $u\in\mathcal T$, put $d=fu$,
and define the polynomial
\[
 P(T)=e^2w(d)T^2-d^2w(e)T-d^2e.
\]
Its leading coefficient is positive, so $P$ is nonzero
and has degree two.
For $t\in\mathcal T$, we have $t>e$, and therefore
$w(et)=tw(e)+e$ by Lemma \ref{lemmawap}. It follows that
\[
 P(t)=H(et,d).
\]
The content of $P$ is
\begin{align*}
 \operatorname{cont}(P)
 &=\gcd(e^2w(d),d^2w(e),d^2e)
 \\
 &=\gcd(e^2w(d),d^2)
 \\
 &=\gcd(e^2,d^2)
 \\
 &=\gcd(e,d)^2
 =\gcd(e,f)^2.
\end{align*}
By Lemma~\ref{gcdmwm}, applied with $n=e$ and $n=d$,
we have $\gcd(w(e),e)=\gcd(w(d),d)=1$.
These identities justify the second and third equalities,
respectively. The last equality uses the fact that $u$ is
prime and $u>e$, so that $\gcd(e,u)=1$.

 For each real $v$ with $Y^{10} < v \leq Y^{11}$, we apply the upper bound
 in Lemma \ref{lemboundw} with $n=d$ and $n=e$. 
 Using $e \leq Y^{2}$, $d \leq Y^{13}$, and $v \leq Y^{11}$, we obtain
\begin{align*}
 |P(v)| & \leq e^2w(d)v^2+d^2w(e)v+d^2e
 \\
 &\leq\frac{15}{\log2}Y^{39}\log Y+Y^{28}
 \leq Y^{40}
\end{align*}
for all sufficiently large $Y$, uniformly in $e,f,u$.
Thus the polynomial averaging estimate \eqref{llmany} applies
with
$$
 \alpha=10,\quad\delta=11,\quad\beta=\tfrac12,
 \quad K=14,\quad M=40,
$$
 and degree bound two. For $t\in\mathcal T$, the preceding nonvanishing argument
shows that $P(t)=0$ if and only if $et=fu$.
 This gives us that 
\begin{align}
\begin{split}
 \sum_{\substack{t\in\mathcal T\\et\ne fu}}
 \frac{\mathcal B_{\frac{1}{2},14}(|H(et,fu)|)}t
 &\ll\mathcal B_{\frac{1}{2},14}(\gcd(e,f)^2)
 \\
 &=\mathcal B_{\frac{1}{2},14}(\gcd(e,f)).
\end{split}\label{intermediatet}
\end{align}
The last equality holds because an integer and its square
have exactly the same prime divisors.

The implied constant in \eqref{intermediatet} is independent of
$e$, $f$, and $u$. Indeed, the parameters in \eqref{llmany}
have the fixed values specified above, and the bound
$|P(v)|\leq Y^{40}$ holds uniformly for all admissible choices
of these factors. Thus \eqref{intermediatet} holds
simultaneously for all such choices once $Y$ exceeds one
absolute threshold. We may therefore multiply by
$\frac{1}{efu}$ and sum over $e,f\in\mathcal E$ and
$u\in\mathcal T$, with the dependence on $e$ and $f$
retained in the explicit factor
$\mathcal B_{\frac{1}{2},14}(\gcd(e,f))$.

 By the uniqueness asserted for the factorization in
\eqref{cequalset}, each $c\in\mathcal C$ determines exactly
one pair $(e,t)\in\mathcal E\times\mathcal T$ satisfying
$c=et$, and each $d\in\mathcal C$ determines exactly one
pair $(f,u)\in\mathcal E\times\mathcal T$ satisfying $d=fu$.
Conversely, each such choice of factors determines
$c,d\in\mathcal C$.
Since $c\ne d$ if and only if $et\ne fu$, and $cd=eftu$,
we obtain
\begin{align*}
 & \sum_{\substack{c,d\in\mathcal C\\c\ne d}}
 \frac{ \mathcal B_{\frac{1}{2},14}(|H(c,d)|) }{cd}
 \\
 &\qquad=
 \sum_{e,f\in\mathcal E}\frac1{ef}
 \sum_{u\in\mathcal T}\frac1u
 \sum_{\substack{t\in\mathcal T\\et\ne fu}}
 \frac{\mathcal B_{\frac{1}{2},14}(|H(et,fu)|)}t
 \\
 &\qquad\ll
 \left(\sum_{u\in\mathcal T}\frac1u\right)
 \sum_{e,f\in\mathcal E}
 \frac{\mathcal B_{\frac{1}{2},14}(\gcd(e,f))}{ef}
 \\
 &\qquad\ll(\log Y)^2.
\end{align*}
 Here \eqref{mainprimemass} bounds the sum over $u$. For the remaining sum, all terms are nonnegative, so we
include all positive integers $e,f\leq Y^2$ and apply the
harmonic gcd averaging lemma with
$\beta = \frac{1}{2}$, $K=14$, and $Z=Y^2$.
Combining this estimate with \eqref{mainconditionalG} gives
\begin{equation}\label{maingcdmoment}
 \sum_{\substack{c,d\in\mathcal C\\c\ne d}}
 \sum_{r,s\in\mathcal R}
 \frac{G(cr,ds)^{\frac{3}{2}}}{cdrs}
 \ll(\log Y)^2.
\end{equation}

 For every $p\geq2$, we have 
$
 \big(1+\frac2p \big)^3
 \leq1+\frac{14}{p} $ (as may be seen by expanding the left-hand side). 
Multiplying over the prime divisors of a positive integer $N$
gives
\begin{equation}\label{mainBcube}
 \mathcal B_{1,2}(N)^3\leq\mathcal B_{1,14}(N).
\end{equation}

Fix $c,d\in\mathcal C$ and $s\in\mathcal R$, allowing $c=d$,
and put $U=w(c)$ and $W=w(d)$. Consider
\begin{equation}\label{defineFT}
 F(T)=(UT+c)(cT-ds).
\end{equation}
This is a polynomial of degree two with positive leading
coefficient $Uc$.
Its coefficients are $Uc$, $c^2-Uds$, and $-cds$.
Since $U=w(c)$, Lemma~\ref{gcdmwm}, applied with
$n=c$, gives $\gcd(U,c)=1$. Consequently, the equalities among 
\[
 \gcd(U,c^2-Uds)=\gcd(U,c^2)=1 
\]
 hold. 
It follows that
\begin{align*}
 \operatorname{cont}(F)
 &=\gcd(Uc,c^2-Uds,cds)
 \\
 &=\gcd(c,c^2-Uds,cds)
 \\
 &=\gcd(c,Uds)
 \\
 &=\gcd(c,ds)
 =\gcd(c,d).
\end{align*}
In the second equality, the factor $U$ can be removed because
it is relatively prime to $c^2-Uds$. The third equality uses
$c\mid cds$ and $c^2-Uds\equiv-Uds\pmod c$.
The fourth uses $\gcd(c,U)=1$, and the last uses $\gcd(c,s)=1$,
since $s$ is prime and $s>c$.

Applying the upper bound in Lemma~\ref{lemboundw} with
$n=c$ and $n=d$, and using $U=w(c)$, $W=w(d)$, and
$c,d\leq Y^{13}$, we obtain
\[
 U,W\leq\frac{13}{\log2}Y^{13}\log Y\leq Y^{14}
\]
for all sufficiently large $Y$, uniformly for $c,d\in\mathcal C$.
 For every real $v$ with $Y^{40}<v\leq Y^{41}$, we therefore have
\[
 0<Uv+c\leq2Y^{55},
 \qquad
 |cv-ds|\leq Y^{54}.
\]
The second inequality holds because both $cv$ and $ds$ lie
between zero and $Y^{54}$.
Thus
\[
 |F(v)|\leq2Y^{109}\leq Y^{110}.
\]
The polynomial averaging estimate \eqref{llmany} consequently
applies with $\alpha=40$, $\delta=41$, $\beta=1$, $K=14$,
degree bound two, and $M=110$, uniformly in $c,d,s$.

 For $r\in\mathcal R$, the factor $Ur+c$ is positive, so, 
 from the definition of $F$ in \eqref{defineFT}, 
 we have that $
 F(r)$ is nonvanishing if and only if $cr$ and $ ds$ are distinct. 
 Moreover, whenever $cr\ne ds$, we also have that 
$
 D(cr,ds)=(Ws+d)|F(r)|$.
 Using Lemma~\ref{lemcalBineq} and
\eqref{llmany} together, this gives us that 
\begin{align*}
 \sum_{\substack{r\in\mathcal R\\cr\ne ds}}
 \frac{\mathcal B_{1,14}(D(cr,ds))}{r}
 &\leq
 \mathcal B_{1,14}(Ws+d)
 \sum_{\substack{r\in\mathcal R\\F(r)\ne0}}
 \frac{\mathcal B_{1,14}(|F(r)|)}r
 \\
 &\ll
 \mathcal B_{1,14}(Ws+d)
 \mathcal B_{1,14}(\gcd(c,d)).
\end{align*}
 Since $W=w(d)$, Lemma~\ref{gcdmwm}, applied with
$n=d$, gives us that $\gcd(W,d)=1$, 
 so that the polynomial
$WT+d$ is primitive. Moreover,
\[
 0<Wv+d\leq2Y^{55}\leq Y^{56}
 \qquad(Y^{40}<v\leq Y^{41}).
\]
 Applying \eqref{llmany} to this polynomial with degree bound one and $M=56$, with the other parameters unchanged, gives
$$
 \sum_{s\in\mathcal R}\frac{\mathcal B_{1,14}(Ws+d)}s\ll1.
$$
Multiplying the preceding bound for the sum over $r$ by $\frac{1}{s}$
and summing over $s\in\mathcal R$ consequently yields
\begin{equation}\label{mainconditionalB}
 \sum_{\substack{r,s\in\mathcal R\\cr\ne ds}}
 \frac{\mathcal B_{1,14}(D(cr,ds))}{rs}
 \ll\mathcal B_{1,14}(\gcd(c,d)).
\end{equation}

The implied constant and the lower bound on $Y$ in
\eqref{mainconditionalB} are independent of $c$ and $d$.
Both applications of \eqref{llmany} use fixed numerical
parameters and fixed degree bounds, and their real-interval
estimates hold uniformly over the remaining factors.
We may therefore sum \eqref{mainconditionalB} with weights
$\frac{1}{cd}$ over $c,d\in\mathcal C$.

Applying the uniqueness assertion accompanying
\eqref{cequalset} to both $c,d\in\mathcal C$, we obtain
unique factors $e,f\in\mathcal E$ and $t,u\in\mathcal T$
such that $c=et$ and $d=fu$.
Both $t$ and $u$ exceed both $e$ and $f$.
If $t\ne u$, then $\gcd(t,fu)=\gcd(u,et)=1$, so
$\gcd(et,fu)=\gcd(e,f)$.
If $t=u$, then $\gcd(et,fu)=t\gcd(e,f)$.
Hence
\[
 \gcd(et,fu)=
 \begin{cases}
 \gcd(e,f),&t\ne u,\\
 t\gcd(e,f),&t=u.
 \end{cases}
\]
 In the $t=u$ case, applying Lemma~\ref{lemcalBineq}
with $\beta=1$, $K=14$, $M=t$, and $N=\gcd(e,f)$,
and using $\mathcal{B}_{1,14}(t)=1 + \frac{14}{t}$ and $t\geq2$,
we obtain
\[
 \mathcal B_{1,14}(t\gcd(e,f))
 \leq\left(1+\frac{14}{t}\right)
 \mathcal B_{1,14}(\gcd(e,f))
 \leq8\mathcal B_{1,14}(\gcd(e,f)).
\]
The same bound therefore holds in both cases, and we obtain
\begin{align*}
 \sum_{c,d\in\mathcal C}
 \frac{\mathcal B_{1,14}(\gcd(c,d))}{cd}
 &=
 \sum_{e,f\in\mathcal E}\frac1{ef}
 \sum_{t,u\in\mathcal T}
 \frac{\mathcal B_{1,14}(\gcd(et,fu))}{tu}
 \\
 &\leq
 8\left(\sum_{t\in\mathcal T}\frac1t\right)^2
 \sum_{e,f\in\mathcal E}
 \frac{\mathcal B_{1,14}(\gcd(e,f))}{ef}
 \\
 &\ll(\log Y)^2.
\end{align*}
The final inequality uses \eqref{mainprimemass} and the
harmonic gcd averaging lemma with $\beta=1$, $K=14$, and
$Z=Y^2$, after including all positive integers $e,f\leq Y^2$.
Summing \eqref{mainconditionalB} with weights $\frac{1}{cd}$ and
using \eqref{mainBcube} now gives
\begin{equation}\label{mainBmoment}
 \sum_{\substack{a,b\in\mathcal A\\a\ne b}}
 \frac{\mathcal B_{1,2}(D(a,b))^3}{ab}
 \ll(\log Y)^2.
\end{equation}
By the uniqueness asserted for the factorization $a=cr$
in \eqref{aequalscr}, applied to both $a$ and $b$,
each ordered pair $(a,b)\in\mathcal A^2$ determines exactly
one choice of $c,d\in\mathcal C$ and $r,s\in\mathcal R$
such that $a=cr$ and $b=ds$.
Under these substitutions, the condition $cr\ne ds$
is equivalent to $a\ne b$.

 Applying the formulation of H\"older's inequality 
 indicated in \eqref{generalHolder} with exponents $\frac{3}{2}$ and $3$
to the finite sum over $c,d\in\mathcal C$, $c\ne d$,
and $r,s\in\mathcal R$, we obtain
\begin{multline*}
 \sum_{\substack{c,d\in\mathcal C\\c\ne d}}
 \sum_{r,s\in\mathcal R}
 \frac{G(cr,ds)\mathcal B_{1,2}(D(cr,ds))}{cdrs}
 \leq \\ 
 \left(
 \sum_{\substack{c,d\in\mathcal C\\c\ne d}}
 \sum_{r,s\in\mathcal R}
 \frac{G(cr,ds)^{\frac{3}{2}}}{cdrs}
 \right)^{\frac{2}{3}}
 \left(
 \sum_{\substack{c,d\in\mathcal C\\c\ne d}}
 \sum_{r,s\in\mathcal R}
 \frac{\mathcal B_{1,2}(D(cr,ds))^3}{cdrs}
 \right)^{\frac{1}{3}}.  
\end{multline*}
 The first factor is bounded according to \eqref{maingcdmoment}.
 For the latter factor, we find that $c\ne d \Longrightarrow cr\ne ds$ by
 virtue of \eqref{mainuniqueab}, so that, with the use of \eqref{mainBmoment}, we obtain that 
\begin{align*}
 \sum_{\substack{c, d \in \mathcal{C} \\ c \neq d}} \sum_{r, s \in \mathcal{R}}
 \frac{\mathcal{B}_{1, 2}\big( D(cr, ds) \big)^{3}}{cdrs} 
 & \leq \sum_{c, d \in \mathcal{C}} \sum_{\substack{r, s \in \mathcal{R} \\ cr \neq ds}}
 \frac{\mathcal{B}_{1, 2}\big( D(cr, ds) \big)^3}{cdrs} \\ 
 & = \sum_{\substack{a, b \in \mathcal{A} \\ a \neq b}} \frac{\mathcal{B}_{1, 2}\big( D(a, b) \big)^{3}}{ab} \\ 
 & \ll \left( \log Y \right)^{2}. 
\end{align*} 
 The foregoing arguments combine to give us that 
\begin{equation}\label{mainweightedmoment}
 \sum_{\substack{c,d\in\mathcal C\\c\ne d}}
 \sum_{r,s\in\mathcal R}
 \frac{G(cr,ds)\mathcal B_{1,2}(D(cr,ds))}{cdrs} \ll(\log Y)^2.   
\end{equation}

 For each $a\in\mathcal A$, we have $2\parallel a$ by \eqref{maintwodivides}. Applying both bounds in Lemma~\ref{lemboundw} with 
 $n=a$, we obtain
\begin{equation}\label{mainwsize}
 \frac a2\leq w(a)
 \leq\frac{a\log a}{\log2}
 \leq\frac{54}{\log2}Y^{54}\log Y
 \leq Y^{55}=x^{0.055}
\end{equation}
 for all sufficiently large $Y$, uniformly in $a$. Also, $a\leq Y^{54}=x^{0.054}$ by \eqref{mainfactorranges}.
For distinct $a,b\in\mathcal A$, the bounds in
\eqref{mainwsize}, applied to both $a$ and $b$, together with
$a,b\leq x^{0.054}$, show that Lemma~\ref{decimallemma}
applies with $A=w(a)$ and $B=w(b)$.

Fix $c,d\in\mathcal C$ with $c\ne d$ and $r,s\in\mathcal R$,
and put $a=cr$ and $b=ds$. For primes in the defining
intervals for $ap,bq\in\mathcal N_x$, we have $p\nmid a$
and $q\nmid b$, so Lemma \ref{lemmawap} gives the biconditional equivalence such that 
 $
 w(ap)=w(bq)
 \Longleftrightarrow 
 w(a)p-w(b)q=b-a$. 
 If $\gcd(w(a),w(b))$ does not divide $b-a$, there are no
such pairs. If it does divide $b-a$, then
$
 G(a,b)=\gcd(w(a),w(b))$. In both cases, Lemma~\ref{decimallemma} bounds the number
of ordered prime pairs by
\begin{equation}\label{applicationdecimal}
 C\left(
 \frac{xG(a,b)}{ab(\log x)^2}\mathcal B_{1,2}(D(a,b))
 +x^{0.03}
 \right),
\end{equation}
 with one absolute constant $C$.

 Now, for fixed $c, d \in \mathcal{C}$ and $r, s \in \mathcal{R}$,
 let $M_{x}(c, d, r, s)$ denote the cardinality of the set of all pairs
 $(p, q)$ such that $p$ and $q$ are prime, $\frac{x}{2cr} < p \leq \frac{x}{cr}$,
 $\frac{x}{2ds} < q \leq \frac{x}{ds}$, and $w(crp) = w(dsq)$.
 Also, let $E_{\neq}(x)$ denote the number of ordered pairs
$(n_1,n_2)\in\mathcal N_x^2$ satisfying $w(n_1)=w(n_2)$
for which the unique factorizations $n_1=crp$ and $n_2=dsq$
have $c\ne d$.
 The uniqueness established in \eqref{aequalscr} and \eqref{nequalsap} then
 gives us the counting identity such that
\begin{equation}\label{Equadruple}
 E_{\neq}(x) = \sum_{\substack{c, d \in \mathcal{C} \\ c \neq d}} \sum_{r, s \in \mathcal{R}} M_x(c, d, r, s).
\end{equation}
For $c,d\in\mathcal C$ with $c\ne d$ and $r,s\in\mathcal R$,
we have $cr\ne ds$ by \eqref{mainuniqueab}.
The application of Lemma~\ref{decimallemma} recorded in
\eqref{applicationdecimal}, with $a=cr$ and $b=ds$, gives
\begin{equation}\label{boundMxcdrs}
 M_{x}(c, d, r, s) \leq
 C \left( \frac{x}{\left( \log x \right)^2 } \frac{G(cr, ds) \mathcal{B}_{1, 2}\big( D(cr, ds) \big) }{ c d r s }
 + x^{0.03} \right).
\end{equation}

 With regard to the index sets associated with the quadruple sum in \eqref{Equadruple},
 the number of corresponding $4$-tuples
 is at most
 $| \mathcal{C} |^{2} |\mathcal{R}|^{2} = |\mathcal{A}|^2$,
 and this is because the unique factorization property
 in \eqref{aequalscr} gives us that $|\mathcal{A}| = |\mathcal{C}| |\mathcal{R}|$.
 This and \eqref{Equadruple} and \eqref{boundMxcdrs} together
 give us that
\begin{equation}\label{boundEneq}
 E_{\neq}(x) \leq
 \frac{Cx}{\left( \log x \right)^{2} }
 \sum_{\substack{ c, d \in \mathcal{C} \\ c \neq d }} \sum_{r, s \in \mathcal{R}}
 \frac{G(cr, ds) \mathcal{B}_{1, 2}\big( D(cr, ds) \big)}{c d r s} + C |\mathcal{A}|^{2} x^{0.03}.
\end{equation}

Since $\mathcal A$ consists of positive integers at most $Y^{54}$,
we have $|\mathcal A|\leq Y^{54}$.
Applying \eqref{mainweightedmoment} to the right-hand side of
\eqref{boundEneq}, we obtain
\begin{align}
 & \frac{Cx}{(\log x)^2}
 \sum_{\substack{c,d\in\mathcal C\\c\ne d}}
 \sum_{r,s\in\mathcal R}
 \frac{G(cr,ds)\mathcal B_{1,2}(D(cr,ds))}{cdrs}
 +C|\mathcal A|^2x^{0.03}
 \nonumber \\
 &\ll\frac{x(\log Y)^2}{(\log x)^2}
 +Y^{108}x^{0.03}
 \nonumber \\
 &\ll x+Y^{138} \nonumber \\
 & \ll x. \nonumber
\end{align}
Here $x=Y^{1000}$, so $x^{0.03}=Y^{30}$ and
$Y^{138}=o(x)$.

Consequently, the relation
\begin{equation}\label{maindifferentcount}
 E_{\neq}(x) \ll x
\end{equation}
 holds.

 Fix $c\in\mathcal C$ and distinct $r,s\in\mathcal R$, and let $p$ and $q$ be primes such that $\frac{x}{2cr} < p \leq \frac{x}{cr}$ and 
 $\frac{x}{2cs} < q \leq \frac{x}{cs}$. Since $cr,cs\leq Y^{54}$ and $p,q>\tfrac{1}{2}Y^{946}>Y^{54}$ for sufficiently large $x$, we have 
 $p\nmid cr$ and $q\nmid cs$. Using the identities
$w(crp)=p\,w(cr)+cr$ and $w(csq)=q\,w(cs)+cs$, we find that
the equality $w(crp)=w(csq)$ is equivalent to
\begin{equation}\label{mainsameequation}
 w(cr)p-w(cs)q=c(s-r).
\end{equation}
If this equation has no integer solutions, it contributes
no prime pairs. Otherwise, fix one integer solution
$(p_0,q_0)$. For any other integer solution $(p,q)$,
subtracting the two equations gives
$
 w(cr)(p-p_0)=w(cs)(q-q_0)$. 
The integers
$\frac{w(cr)}{\gcd(w(cr),w(cs))}$ and
 $\frac{w(cs)}{\gcd(w(cr),w(cs))}$ are relatively prime.
It follows that
 $ \frac{w(cs)}{\gcd(w(cr),w(cs))}\mid (p-p_0)$. Thus the possible values of $p$ lie in one residue class
modulo $ \frac{w(cs)}{\gcd(w(cr),w(cs))}$. Now, observe that the interval
$ \frac{x}{2cr}<p\leq\frac{x}{cr}
$ is of length $\frac{x}{2cr}$. Distinct integers in that residue
class differ by at least the modulus, so the interval
contains at most $ \frac{x\gcd(w(cr),w(cs))}{2cr\,w(cs)}+1 $ possible values of $p$. Each determines at most one $q$
through \eqref{mainsameequation}.
 Since $2\parallel cs$, the lower bound in
Lemma~\ref{lemboundw}, applied with $n=cs$, gives
$w(cs)\geq \frac{cs}{2}$. 
 Therefore the number of ordered prime pairs in the required
intervals satisfying \eqref{mainsameequation} is at most
\begin{equation}\label{mainsamepairbound}
 \frac{x\gcd(w(cr),w(cs))}{c^2rs}+1.
\end{equation}
This bound also holds when there are no integer solutions.

 For fixed $c\in\mathcal C$, let $U=w(c)$. Since $c>1$, the definition of $w$ gives us that $U>0$. By Lemma~\ref{gcdmwm}, applied 
 with $n=c$, we have that $\gcd(U,c)=1$. Moreover, the equality $ \gcd(w(cr),w(cs))=\gcd(Ur+c,Us+c) $
 holds. Since $r,s>Y^{40}$, Lemma~\ref{lemlineargcd}, applied with
$R=Y^{41}$, gives
\begin{align*}
 \sum_{\substack{r,s\in\mathcal R\\r\ne s}}
 \frac{\gcd(w(cr),w(cs))}{rs}
 &\leq
 Y^{-80}
 \sum_{\substack{1\leq m,n\leq Y^{41}\\m\ne n}}
 \gcd(Um+c,Un+c)
 \\
 &\ll Y^{-80}Y^{82}\log(Y^{41})
 \ll Y^2\log Y.
\end{align*}
The constant is independent of $c$ and $U=w(c)$.
Moreover, by \eqref{mainfactorranges}, the relations among 
\begin{align*}
 \sum_{c\in\mathcal C}\frac1{c^2}
 & \leq \sum_{n>Y^{11}}\frac1{n^2} \\ 
 & \leq\int_{\lfloor Y^{11}\rfloor}^{\infty}\frac{dv}{v^2} \\ 
 & = \frac1{\lfloor Y^{11}\rfloor} \\ 
 & \ll Y^{-11} 
\end{align*}
 hold, where $n$ ranges over positive integers. We also find that $|\mathcal C|\leq Y^{13}$ and that $|\mathcal R|\leq Y^{41}$. 
 Summing \eqref{mainsamepairbound} over $c\in\mathcal C$
and distinct $r,s\in\mathcal R$ consequently shows that
the number of ordered equal-value pairs with $c=d$
and $r\ne s$ is at most
 $$ x\sum_{c\in\mathcal C}\frac1{c^2}
 \sum_{\substack{r,s\in\mathcal R\\r\ne s}}
 \frac{\gcd(w(cr),w(cs))}{rs}
 +|\mathcal C|\,|\mathcal R|^2 
 \ll xY^{-11}Y^2\log Y+Y^{13+82}. $$
 As a consequence, we have that 
\begin{multline*}
 \text{($\#$ of ordered equal-value pairs s.t.\ $c = d$ and $r \neq s$)} 
 \ll \\ 
 xY^{-9}\log Y+Y^{95}. 
\end{multline*}
 From this and the equality $x=Y^{1000}$ equivalent to \eqref{x1000th}, 
 we find that 
\begin{equation}\label{mainsamecount}
 \text{($\#$ of ordered equal-value pairs s.t.\ $c = d$ and $r \neq s$)} 
 =o(x). 
\end{equation}

 Finally, if $c=d$ and $r=s$, then $a=b$. The equality $w(ap)=w(aq)$ then yields $w(a)(p-q)=0$. Since $w(a)>0$, it follows that $p = 
 q$. These are precisely the pairs $(n,n)$ with $n\in\mathcal N_x$.
There are 
\begin{equation}\label{suchpairs}
 |\mathcal N_x|\leq x
\end{equation}
 such pairs.

 Define $ E(x)= \left| \{(n_1,n_2)\in\mathcal N_x^2:w(n_1)=w(n_2)\} \right|$. For each ordered pair $(n_1,n_2)\in\mathcal N_x^2$ with 
 $w(n_1)=w(n_2)$, there exist elements $c$ and $d$ in $\mathcal C$, elements $r$ and $s$ in $\mathcal R$, and primes $p$ and 
 $q$ such that $n_1 = c r p$, $n_2=dsq$, $\frac{x}{2cr}<p\leq\frac{x}{cr}$, and $\frac{x}{2ds}<q\leq\frac{x}{ds}$. By \eqref{mainfactorranges}, 
 we have that $cr,ds\leq Y^{54}$. Since $x=Y^{1000}$ and $Y\geq2$, it follows that  $p,  q> 
 \frac12Y^{946}>Y^{54}$. Thus $p$ and $q$ are the largest prime divisors of
$n_1$ and $n_2$, respectively.  Consequently, the primes $p$ and $q$ and the integers $cr = \frac{n_1}{p}$ 
 and $ds = \frac{n_2}{q}$ are uniquely determined by $(n_1,n_2)$.
 By \eqref{mainuniqueab}, the pairs $(c,r)$ and $(d,s)$
are then uniquely determined. Therefore, for each such ordered pair $(n_1,n_2)$,
exactly one of the following conditions holds:
$c\ne d$; $c=d$ and $r\ne s$; or $c=d$ and $r=s$. 

 We have thus partitioned the pairs counted by $E(x)$ into three disjoint classes, writing 
\begin{multline*}
 E(x) = E_{\neq}(x)
 + \\ \text{number of ordered equal-value pairs with $c=d$ and $r\ne s$}
 +|\mathcal N_x|.
\end{multline*}
 By \eqref{mainsamecount}, the preceding equality gives us that $$ E(x) = E_{\neq}(x) + o(x) + |\mathcal{N}_{x}|. $$ From the relation 
 $E_{\neq}(x) \ll x$ highlighted in \eqref{maindifferentcount}
 and from the bound for $|\mathcal{N}_{x}|$ 
 highlighted in \eqref{suchpairs}, we thus obtain that 
\begin{equation}\label{maincollisionbound}
 E(x)\ll x.
\end{equation}

 For each nonnegative integer $v$, define $ R_x(v)=|\{n\in\mathcal N_x:w(n)=v\}|$. Only finitely many $R_x(v)$ are nonzero, with 
 $ \sum_{v\geq0}R_x(v)=|\mathcal N_x|$ and $ \sum_{v\geq0}R_x(v)^2=E(x)$. The second equality holds because $R_x(v)^2$ counts the 
 ordered pairs of elements of $\mathcal N_x$ that both have value $v$ under $w$. 

 Let $v_1$, $v_2$, $\ldots$, $v_k$ be the distinct values attained by $w$ on $\mathcal{N}_{x}$. Applying the Cauchy--Schwarz inequality to 
 the vectors $\big( R_x(v_1), R_x(v_2), \ldots, R_x(v_k) \big)$ and $(1, 1, \ldots, 1)$, we obtain that 
\begin{align*}
 |\mathcal N_x|^2
 &=\left(\sum_{\substack{v\geq0\\R_x(v)>0}}R_x(v)\right)^2
 \\
 &\leq
 |\{v\geq0:R_x(v)>0\}|
 \sum_{v\geq0}R_x(v)^2
 \\
 &=|\{w(n):n\in\mathcal N_x\}|\,E(x).
\end{align*}
 By \eqref{definecalNx} and \eqref{maininterval} together, each element in $\mathcal N_x$ is a positive integer that does not exceed 
 $x$. Also, \eqref{mainmass} implies that $\mathcal N_x$ is nonempty for all sufficiently large $x$, so $E(x)\geq|\mathcal N_x| > 
 0$. Hence $$ V(x)\geq|\{w(n):n\in\mathcal N_x\}| \geq\frac{|\mathcal N_x|^2}{E(x)}. $$ By \eqref{mainmass} and 
 \eqref{maincollisionbound}, there exist absolute constants $c_1, 
 c_2>0$ such that, for all sufficiently large real $x$, each of the relations among $ |\mathcal N_x|\geq c_1x$ 
 and $ E(x)\leq c_2x$ holds. It follows that $ V(x)\geq\frac{c_1^2}{c_2}\,x$. 
Taking $c_0 = \frac{c_1^2}{c_2}$ and choosing one absolute $x_0$
large enough for all the preceding estimates proves the
asserted bound for every real $x\geq x_0$.
Dividing by $x$ and taking the lower limit gives us that \eqref{strongest} holds, 
 as desired. 
\end{proof}

\section{Conclusion}
 Since we have proved, as in Theorem \ref{maintheorem}, that $\liminf_{x\to\infty}\frac{V(x)}x$ is strictly positive, this raises the question as 
 to how $\liminf_{x\to\infty}\frac{V(x)}x$ can be evaluated explicitly or numerically estimated in a meaningful way (noting that it is 
 immediate that $\liminf_{x\to\infty}\frac{V(x)}x \leq 1$). Omitting details, our construction above, as in our proof of 
 Theorem \ref{maintheorem}, can be used to show that $$ \liminf_{x\to\infty}\frac{V(x)}x > 3 \times 10^{-20}. $$ We very much 
 encourage future research endeavors based on the problems of bounding, numerically estimating, and evaluating
 $\liminf_{x\to\infty}\frac{V(x)}x$. 
 This also raises a further open question: Does the limit 
 $\lim_{x\to\infty}\frac{V(x)}x$ exist? 

\begin{figure}
\begin{center}
\includegraphics[scale=0.25]{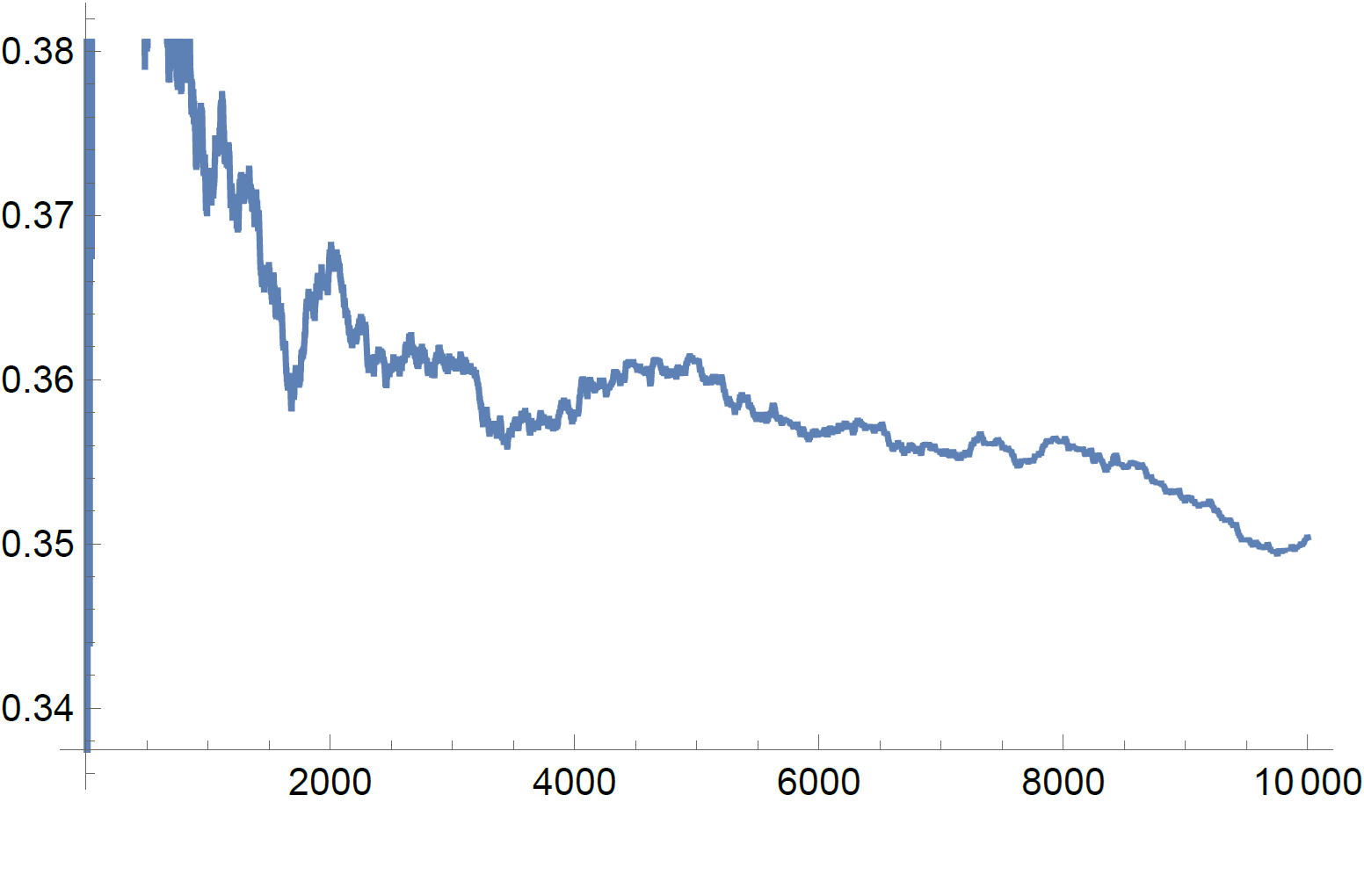}
\caption{\label{graph10000p} A plot of the values of $\frac{V(x)}{x}$ for integers $x \in [1, 10{\,}000]$, 
 generated using the {\tt ListLinePlot} command in Mathematica.}
\end{center}
\end{figure}

 The behavior of $\frac{V(x)}{x}$ for integers $x \in [1, 10{\,}000]$ is illustrated in the plot indicated in Figure \ref{graph10000p}. 

 \subsection*{AI disclosure}
 The following AI disclosure is based on AI disclosures used and written by the author previously (with numerous parts taken verbatim 
 from previous AI disclosures used and written by the author). 

 Our solution above to Erd\H{o}s's problem is, as indicated above, based on our extensive interactions with GPT-6 Astra Pro, and this is 
 described in more detail below. All AI-generated suggestions were substantially revised, corrected, and independently verified by the 
 author, who assumes full responsibility for the mathematical content. The author's use of AI is fully in accordance with editorial standards 
 on the responsible, ethical, and transparent use of AI. 

 The author inputted a prompt into GPT-6 Astra Pro containing an equivalent version of the problem from Guy's text solved above (with an 
 explicit attribution to Guy's text on unsolved problems in number theory), concluding with an explicit request that 
 GPT-6 Astra Pro solve the given problem. After {\tt 16m 16s} of thinking, GPT-6 Astra Pro provided a response beginning with
 the following, and provided a very rudimentary outline of a potential proof: 
\begin{quote}
``{\bf I obtain a negative answer.} The argument below gives the stronger conclusion
 $\liminf_{x \to \infty} \frac{V(x)}{x} > 0.$
 In particular, $V(x) \neq o(x)$.

This is a {\bf proposed proof developed here}. Because it would resolve the published question, it should receive independent specialist review before submission. I have written out the uniform estimates and the exceptional cases on which the argument depends.''
\end{quote}

 After extensive clarifications, corrections, revisions, suggestions, verifications, calculations, etc., through extensive discussions between 
 the author and GPT-6 Astra Pro, and after extensive and careful reviews of such discussions on the part of the author, this provided a 
 basis for the above solution, written by the author, to Erd\H{o}s's problem. 
 
 Specific examples of corrections, revisions, clarifications, constructions, etc., on the part of the author include a modification of a previous 
 version of the bound in \eqref{x0geq8} to construct a corrected proof associated with the $\rho(2) = 2$ case involved in the proof of 
 Lemma \ref{decimallemma}, the use of \eqref{simplifyint} in relation to \eqref{splitloglog}--\eqref{logYpoweroverh}, the use of the bound 
 in \eqref{boundBasel} for quadratic harmonic numbers, the use of the bound in \eqref{classicalskew} for skew harmonic numbers, the use 
 of an upper bound for harmonic numbers in relation to \eqref{floorharmonic}, the many uses of and references to Lemmas 
 \ref{gcdmwm}--\ref{lemcalBineq} throughout the paper, the justification related to the derivation of \eqref{maincollisionbound} as to how 
 the pairs $(c,r)$ and $(d,s)$ are uniquely determined by invoking \eqref{mainuniqueab}, the justification based on the uniqueness of 
 the factorization in \eqref{cequalset} and related to the bounds for $ \sum_{\substack{c,d\in\mathcal C\\ c \neq d}}
 \frac{\mathcal B_{\frac{1}{2},14}(|H(c,d)|)}{cd}$, the relation in \eqref{mainsamecount} and its derivation, the use of and derivation of the 
 inequality $$ \sum_{\substack{c, d \in \mathcal{C} \\ c \neq d}} \sum_{r, s \in \mathcal{R}} \frac{\mathcal{B}_{1, 2}\big( D(cr, d 
 s) \big)^{3}}{cdrs} \leq \sum_{\substack{a, b \in \mathcal{A} \\ a \neq b}} \frac{\mathcal{B}_{1, 2}\big( D(a, b) \big)^{3}}{ab}$$
 in the derivation of \eqref{mainweightedmoment}, the counting identity in \eqref{Equadruple} and the use of the bounds in 
 \eqref{boundMxcdrs} and \eqref{boundEneq} to justify both \eqref{maindifferentcount} and \eqref{maincollisionbound}, a correction of 
 an incorrect statement in the proof of Theorem \ref{maintheorem} concerning how the right-hand side of \eqref{mainconditionalG} 
 relates to \eqref{intermediatet} and to the estimate $$ \sum_{\substack{c,d\in\mathcal C\\ c \neq d}} \frac{ \mathcal B_{\frac{1}{2}, 
 14}(|H(c,d)|) }{cd} \ll(\log Y)^2 $$ established in the proof of Theorem \ref{maintheorem}, the use of the consequence in 
 \eqref{furtherCRT} of the Chinese remainder theorem in relation to \eqref{quotientCRT}, \eqref{lcmCRT}, and \eqref{lambdarhoCRT}, 
 the use of \eqref{harmonice} in relation to the derivation of \eqref{mainmass}, etc. 

 Extensive Mathematica programs for numerically verifying the above construction have been carefully prepared and written and tested by 
 the author 
 and are available upon request. 

\bibliographystyle{plain}
\bibliography{seprefe}

\end{document}